\documentclass[letterpaper, 10 pt, conference]{ieeeconf}  % Comment this line out if you need a4paper

\IEEEoverridecommandlockouts                              % This command is only needed if 
\usepackage{cite}
\usepackage{color}
\usepackage{amsfonts,amssymb}
\usepackage{bm}
\usepackage{amsmath}
\usepackage{multirow}
\usepackage{algorithm, algorithmic}
\usepackage{graphicx}
\usepackage{subfigure}
\usepackage{float}
\usepackage{textcomp}
\usepackage{graphicx} %use graph format
\usepackage{epstopdf}
\usepackage{array}
\usepackage[thmmarks,amsmath]{ntheorem}
\newtheorem{mythm}{Theorem}
\newtheorem{mydef}{Definition}
\newtheorem{assumption}{Assumption}
\newtheorem{lemma}{Lemma}
\newtheorem{remark}{Remark}

\newtheorem{corollary}{Corollary}

\usepackage{array}
\usepackage{enumerate}
\usepackage{booktabs}
\usepackage{algorithm}
\usepackage{algorithmic}
\usepackage{setspace}
\usepackage{float}
\usepackage{cases}
\usepackage{stfloats}

\usepackage{mathrsfs}
\title{\LARGE \bf
A Multi-Player Potential Game Approach for\\
Sensor Network Localization with Noisy Measurements
}

\author{Gehui Xu$^{1}$, Guanpu Chen$^{2}$, Baris Fidan$^{3}$,  Yiguang Hong$^{4}$,  \\
Hongsheng Qi$^{5}$, Thomas Parisini$^{6}$,  and Karl H. Johansson$^{2}$% <-this % stops a space
\thanks{*This work was supported by Swedish Research Council Distinguished Professor Grant 2017-01078, Knut and Alice Wallenberg Foundation Wallenberg Scholar Grant,  Swedish Strategic Research Foundation SUCCESS Grant FUS21-0026, and also supported 
in part by the Digital Futures Scholar-in-Residence Program, by the European Union’s Horizon 2020 research and
innovation programme under grant agreement no. 739551 (KIOS CoE), and
by the Italian Ministry for Research in the framework of the 2017 Program
for Research Projects of National Interest (PRIN), Grant no. 2017YKXYXJ.}% <-this % stops a space
\thanks{$^{1}$Gehui Xu
is with the Key Laboratory of Systems and Control, Academy of Mathematics and Systems Science, Beijing, China, and also  with the Division of Decision and Control Systems, School of	Electrical Engineering and Computer Science, KTH Royal Institute of	Technology, 100 44, Stockholm, Sweden. 
        {\tt\small xghapple@amss.ac.cn}}%
\thanks{$^{2}$ Guanpu Chen, and Karl H. Johansson 
are with the Division of Decision and Control Systems, School of	Electrical Engineering and Computer Science, KTH Royal Institute of	Technology, 100 44, Stockholm, Sweden. 
        {\tt\small  guanpu@kth.se, and kallej@kth.se}}
        \thanks{$^{3}$Baris Fidan is with the Department of Mechanical and Mechatronics Engineering, University of Waterloo, Waterloo, ON N2L 3G1, Canada.
        {\tt\small fidan@uwaterloo.ca}}%
                \thanks{$^{4}$Yiguang Hong is with
		Department of Control Science and Engineering, Tongji University,  Shanghai 201804, China, and is also with Shanghai Research Institute for Intelligent Autonomous Systems, Shanghai 201210 China.
        {\tt\small yghong@iss.ac.cn}}%
                  \thanks{$^{5}$Hongsheng Qi is with Key Laboratory of Systems and Control, Academy of Mathematics and Systems Science, Beijing, China, and also with School of Mathematical Sciences, University of Chinese Academy of Sciences, Beijing, China.
        {\tt\small qihongsh@amss.ac.cn}}
                 \thanks{$^{6}$Thomas Parisini is with  the Department of Electrical and Electronic Engineering,
	Imperial College London, London SW7 2AZ, UK, and also with the Department
	of Engineering and Architecture, University of Trieste, Trieste 34127, Italy.
        {\tt\small t.parisini@imperial.ac.uk}}
}

\begin{document}

\maketitle
\thispagestyle{empty}
\pagestyle{empty}

%%%%%%%%%%%%%%%%%%%%%%%%%%%%%%%%%%%%%%%%%%%%%%%%%%%%%%%%%%%%%%%%%%%%%%%%%%%%%%%%
\begin{abstract}

Sensor network localization (SNL) is a challenging problem due to its inherent non-convexity and the effects of  noise in inter-node ranging measurements and anchor node position.
%a difficult task in which 
% position measurements   are subject to errors. Attributed to these measurement errors and the inherent nonconvexity of SNL itself, obtaining the existence condition of a globally localizable network is challenging. 
We formulate a non-convex SNL problem as a multi-player non-convex potential game and investigate the existence and uniqueness of a Nash equilibrium (NE) in both the ideal setting without measurement noise and the practical setting with measurement noise.
%noiseless cases and noise scenarios.  
We first show that the NE exists and is unique in the noiseless case, and corresponds to the precise network localization. Then, we study the SNL for the case
%consider a noisy localization case 
with errors affecting the anchor node position and the inter-node
%measurements and 
distance measurements.  Specifically, we establish that in case these errors are sufficiently small, the NE exists and is unique. It is shown that the NE is an approximate solution to the SNL problem, and that 
%Moreover, if the measurement errors are constrained by a small bound, 
the position errors can be quantified accordingly. Based on these findings, we apply the results to case studies involving only inter-node  distance measurement errors and only anchor position information inaccuracies. 
%Finally, simulation results are provided to validate the theoretical results  and effectiveness of the proposed approach. 
% Sensor network localization (SNL) problems require determining the physical coordinates of all sensors in a network. This process relies on the global coordinates of anchors and the available measurements between non-anchor and anchor nodes. Attributed to the intrinsic non-convexity, obtaining a globally optimal solution to SNL 
% %network 
% is challenging, as well as implementing corresponding algorithms. In this paper, we formulate a non-convex multi-player potential game for a generic SNL problem to investigate the identification condition of the global Nash equilibrium (NE) therein, where the global NE represents the global solution of SNL. We employ canonical duality theory to transform the non-convex game into a complementary dual problem. Then we develop a conjugation-based algorithm to compute the stationary points of the complementary dual problem. On this basis,
%  %With the guarantee that a global NE is equivalent to a global network’s solution, 
%  we show an identification condition of the global NE: the stationary point of the proposed algorithm satisfies a duality relation. Finally, simulation results are provided to validate the effectiveness of the theoretical results.

\end{abstract}

%%%%%%%%%%%%%%%%%%%%%%%%%%%%%%%%%%%%%%%%%%%%%%%%%%%%%%%%%%%%%%%%%%%%%%%%%%%%%%%%
\section{INTRODUCTION}

% Wireless sensor networks
% (WSNs)
% Wireless sensor networks (WSNs), consisting of a number of sensors spread across a geographical area, find diverse applications in target tracking and detection \cite{marechal2010joint}, environment monitoring \cite{sun2010corrosion}, area exploration \cite{sun2005reliable}, data collection, as well as cooperative robot tasks \cite{jing2018weak}.
% find diverse applications \cite{akyildiz2002wireless,liu2015event} owing to their capabilities in sensing, processing, and communication. 
% %Wireless sensor networks (WSNs) have a wide range of applications \cite{akyildiz2002wireless,liu2015event} due to their capabilities of sensing, processing, and communication,
% These applications include target tracking and detection \cite{marechal2010joint}, environment monitoring \cite{sun2010corrosion}, area exploration \cite{sun2005reliable}, data collection,  as well as cooperative robot tasks \cite{jing2018weak}. 
Accurate information regarding the location of nodes within wireless sensor networks (WSNs) is essential in diverse applications, such as  target tracking and detection \cite{marechal2010joint}, environment monitoring \cite{sun2010corrosion}, area exploration \cite{sun2005reliable}, data collection, as well as cooperative robot tasks \cite{jing2018weak}.
%about the location of nodes is essential in these applications. 
%The availability of accurate information about the location of nodes is essential in many sensor network applications.
A common approach for sensor localization involves utilizing (noisy) ranging information obtained through signal transmission techniques, such as time of arrival, time-difference of arrival, and strength of received radio frequency signals \cite{mao2007wireless}.  Also, there are anchor nodes with known global positions \cite{anderson2010formal}. Then a sensor network localization
(SNL) problem  is defined as follows:  Given the positions of  the anchor nodes of the WSN and the measurable information among each non-anchor node and its neighbors,  find the positions of the rest of non-anchor nodes.

To  describe individual preferences and network interactions, game theory constitutes an effective modeling tool 
%well-suited to model the strategic behavior in wireless sensor network 
\cite{jia2013distributed,bejar2010cooperative}. The Nash equilibrium (NE) is a key solution concept characterizing a profile of stable strategies in which rational players would not choose to deviate from their own strategies~\cite{nash1951non,xu2022efficient}. The SNL problem can be formulated as a game by letting non-anchor nodes be players, their estimated position be strategies, and their positioning accuracy measurements be payoffs.
Specifically,  potential games are well-suited to model 
%the strategic behavior in 
SNL  (see \cite{jia2013distributed,ke2017distributed,xu2024global}).
%Note that the sensors need to consider the positioning accuracy of the whole WSN while ensuring their own positioning accuracy through the given information. 
In particular, the potential game paradigm guarantees an alignment between the individual sensor profits and the network objective by exploiting a  unified potential function.   It is possible to find the  NE corresponding to the global optimum for the whole WSN rather than a local approximation.
% To better describe a WSN and each sensor's possible and ideal localization actions,
% game theory  is found useful in modeling  WSNs and SNL problems \cite{jia2013distributed,bejar2010cooperative,chen2023global}. The Nash equilibrium (NE) is a prominent concept in game theory, which characterizes a  profile of stable strategies  where rational sensor nodes would not choose to deviate from their location strategies \cite{nash1951non,xu2022efficient}. Particularly, 
% %potential game is one of the most important noncooperative game models, 
% potential game is well-suited to model the strategic
% behavior in SNL problems \cite{jia2013distributed,ke2017distributed}.
% Note that the sensors need to consider the positioning accuracy of the whole WSN while ensuring their own positioning accuracy through the given information. The potential game framework can
%  guarantee such an alignment between the individual sensor’s profit and the global network's objective by characterizing a global unified potential function. 
 %Moreover, the global NE, instead of local NE, is equal to the network’s unique localization when the grounded graph of sensor network is generically globally rigid and the measurement information is exact.
 
% It is worth noting that non-convexity is an intrinsic challenge of SNL problems
%Due to the inherent non-convexity of the SNL problem,
It is worth noting, though, that   the effects of the
noises in inter-node ranging measurements and the inaccuracies
in anchor node position information
%inexact measurements
%non-convexity 
are an intrinsic challenge of SNL problems, which, unfortunately, cannot be readily avoided by potential games or other modeling approaches. 
%In practice, the measurement information may never be exact.
The inter-node distance measurements between sensor nodes are often subject to inevitable uncertainties like transmission interference or time delay, leading to  measurement errors \cite{kannan2010analysis}.
Besides, many existing studies assume precise anchor positions to estimate the location of the rest of the sensor nodes.  However, in many scenarios, anchor positions may not be accurately known. This is often attributed to the reliance on 
%the anchor positions are usually obtained through GPS  or other positioning systems, thereby causing position information inaccuracies.
the Global Positioning System (GPS)  or other positioning systems for determining anchor positions, which may cause estimation errors \cite{lui2008semi}. 

A fundamental problem in such noisy SNL problems  is whether a given sensor network can be uniquely localized.
The uniqueness of localizability is usually revealed using graph rigidity theory \cite{anderson2008rigid}. When the measurement information is exact, the generically global rigidity of the grounded graph guarantees the uniqueness of NE, corresponding to the precise localization of the network. However, in the presence of noise, there may exist multiple sets of non-congruent sensor localizations that satisfy the provided inaccurate measurements.
%this problem could be ill-posed as there may exist
% if the actual geometries of some parts of the sensor network
% are sensitive to distance measurement errors, it may cause
% erroneous local geometric realizations in those parts of the
% network.
The presence of noise may perturb the rigidity of the graph, resulting in flip ambiguities or even no solution to the problem. Such a localization problem with measurement errors
%inexact or uncertain information
has been considered in some existing works. \cite{lui2008semi} studied an uncertain SNL problem with Gaussian distributed disturbances and used the semi-definite programming relaxation technique to solve the formulated maximum likelihood estimation problem. Then, \cite{naddafzadeh2013second} considered a similar formulation and employed second-order cone programming to seek a robust solution with complexity reduction. Also, \cite{anderson2010formal} investigated an SNL problem with errors in inter-node distance measurements and transformed it into a minimization problem to reveal the relation between errors in positions and errors in distance measurements.

%we will argue that, given the graph theoretic conditions that would guarantee unique localizability in the noiseless case, localization in the noisy case can be posed as a minimization problem, the solution of which has several properties 

%distance measurements will never be exact, and the equations whose solutions deliver sensor positions in the noiseless case in general no longer have a solution. 

%the global NE, is equal to the network’s unique localization, which is when the grounded graph
 
 %is equal to a global optimum of the potential function denoting the network’s precise localization.
 
% ,  it is natural and essential to seek a global NE of the whole sensor network rather than local NE and approximate solutions, since a global NE is equal to a global optimum of the potential function denoting the network’s precise localization.
 %Employing the potential game framework for describing the interaction and schedule between sensors has become  an emerging research topic \cite{jia2013distributed,ke2017distributed}.

The main objective of this paper is to provide a practical solution to the noisy SNL problem. We
 focus on the case
% question of what happens 
 when the inter-node distance measurements and anchor position information are subject to some errors, presumably on account of the measurement process.
Specifically, we formulate the non-convex SNL problem as a potential game and investigate the existence and uniqueness of NE in both the ideal setting with accurate anchor node location information and accurate inter-node distance measurements and the practical setting
with anchor location inaccuracies and distance measurement noise.    We first show that the NE exists and is unique in the noiseless case, corresponding to the precise network localization. Then we study SNL for the case with errors in anchor node position and inter-node distance measurements. We establish that if these errors are sufficiently small, the NE exists and is unique.  It is shown that the NE is an approximate solution to the SNL problem, and that the position errors can be quantified accordingly. 
% , providing
% an approximate solution to the SNL problem. Moreover, if
% the measurement errors are constrained by a small bound,
% the position errors can be quantified accordingly. 
%Based on these findings, 
Then,
we apply the results to case studies involving
only inter-node distance measurement errors and only anchor
position information inaccuracies. 

\section{Problem formulation}

In this section, we first introduce the noisy inter-node range measurement based SNL problem of interest and then formulate it as a potential game.
%range-based noisy 
%with the design scheme for solving the potential game.
%first introduce the sensor network localization problem 
%in \ref{SNL-problem}, 
%and  
%give some preliminary knowledge of potential game in \ref{pre-potential}, finally, we 
%then 

{
\subsection{Noiseless SNL problem}
%in \ref{for-potential}. 
%\subsection{Sensor network localization problem}\label{SNL-problem}
%	In this section, we describe our network model. 
	 %focus on solving the range-based sensor	network localization problem.
 Consider a static sensor network
%	Assume a sensor network
	 in $\mathbb{R}^{n}$  ($ n = 2  $)  composed of 
	  $ M $ anchor nodes
whose positions are known  and $N$ non-anchor sensor nodes whose positions are
unknown ($ M<N $). Let a  graph $ \mathcal{G}=(\mathcal{N}, \mathcal{E}) $ represent the sensing relationships between sensors, where $\mathcal N$ is the sensor node set and $\mathcal E\subseteq\mathcal N \times \mathcal N$ is the set of edges that represent the sensor node pairs whose range measurements are available.
%the edge set between sensors. 
Specifically,  
$ \mathcal{N}= \mathcal{N}_{s}\cup \mathcal{N}_{a}=\{i\}_{i=1}^{\bar{N}}$, where  $\bar{N}=N+M$, $ \mathcal{N}_{s}=\{1,\dots, N\} $ and $ \mathcal{N}_{a}= \{N+1,\dots,N+M \} $ correspond to the sets of non-anchor nodes and anchor nodes, respectively.  Let 
%$ x_{i}^{\star} \in \mathbb{R}^{n} $ for $i\in \mathcal{N}_{s}$ denote the actual position of the $ i $-th 
$ x_{i}^{\star} \in \mathbb{R}^{n} $ denote the actual position of 
sensor node $i\in \mathcal{N}$ with $\boldsymbol{x}^{\star}\triangleq\operatorname{col}\{x_i\}_{i\in \mathcal{N}}$.
% non-anchor
% node, {and	
% $ x_{l}^{\star}=x_{K+m}^{\star} \in \mathbb{R}^{n} $ for $m\in \{1,2,\dots,M \}$
% %for $k\in \mathcal{N}_{a}=\{N+1,N+2,\dots,N+M \}$ 
% denote the actual position of anchor node $l\in \mathcal{N}_{a}$.}
 For a pair of sensor nodes $ i$ and $ j $, their Euclidean distance
% is denoted as $ e_{ik} $.
%Similarly, the distance between non-anchor node  $i $ and non-anchor node  $ j $ 
is denoted by $ d_{ij}^{\star} $. 
Each sensor has the capability
of sensing range measurements from other sensors within a
fixed range $R_s$, and  $\mathcal{E}=\mathcal{E}_{ss}\cup \mathcal{E}_{as}\cup \mathcal{E}_{aa}$ with 
$ \mathcal{E}_{ss}=\{(i,j)\in\mathcal{N}_{s}\times\mathcal{N}_{s}:\|x_{i}^{\star}-x_{j}^{\star}\|\leq R_s,i\neq j\} $  denoting the edge set  between non-anchor nodes,  	$ \mathcal{E}_{as}=\{(i,l)\in\mathcal{N}_{s}\times\mathcal{N}_{a}:\|x_{i}^{\star}-x_{l}^{\star}\|\leq R_s\} $ denoting the edge set between anchor nodes and non-anchor nodes,  $ \mathcal{E}_{aa}=\{(l,m)\in\mathcal{N}_{a}\times\mathcal{N}_{a},l\neq m\} $ denoting the edge set  between anchor nodes.
%Also, suppose that 
%SNL (\ref{SNL}) considers that  the measurements $d_{ij}$ are noise-free and all anchor positions $ x_{l} $, $l\in \mathcal{N}_{a}$ are accurate.
 The 
{ range-based SNL task in the noiseless case is to determine the accurate positions  of all non-anchor sensor nodes $i\in \mathcal{N}_{s}$ when all anchor node positions $ x_{l}^{\star} $, $l\in \mathcal{N}_{a}$ and  measurements $d_{ij}^{\star}$,  $d_{il}^{\star}$ are given: 
%that is, 
}
\begin{align}\label{SNL}
	& \text{find} \quad x_{1},\dots,x_{N}\in \mathbb{R}^{n}\notag\\
	& \text{s.t.} \quad \|x_{i}-x_{j}\|^2=d_{ij}^{\star 2}, \forall (i,j)\in  \mathcal{E}_{ss},\\
	& \quad\quad\;\|x_{i}-x_{l}^{\star}\|^2=d_{il}^{\star 2}, \forall (i,l)\in \mathcal{E}_{as}. \notag \vspace{-0.4cm}
\end{align}
Denote $\boldsymbol{x}^{\dagger}=\operatorname{col}\{x_{1}^{\dagger}, \dots, x_{N}^{\dagger} \} \in \mathbb{R}^{nN} $ as the solution to noiseless SNL problem \eqref{SNL}.
%, which are the actual position vectors of all non-anchor nodes.

%and $\mathcal{N}_{a}^{i}$ as  the  neighbor set of adjacent anchor nodes of the $ i $-th non-anchor node with $(i,k)\in\mathcal{E}_{as}$. 
 %Also, two sensor nodes can communicate with each other if and only if the distance between the two is smaller than a communication range $R_c$. Here we consider $R_s=R_c$.

\subsection{Noisy SNL problem
%Localization with Anchor Position Uncertainty and Distance Uncertainty
}
%In practice, perfect knowledge of anchor positions may not be available. In many scenarios, the anchor positionsare obtained using global positioning system (GPS) and thus subject to estimation errors.
%\noindent\textbf{Anchor position uncertainty}
In practical scenarios, 
%achieving precise knowledge of anchor positions is challenging. The 
the anchor positions are usually obtained through GPS  or other positioning systems, thereby causing position information inaccuracies. For $l \in \mathcal{N}_a$, 
%On this basis,
%$x_l^{\star}$ is the actual position of anchor node $l \in \mathcal{N}_a$ and 
let anchor node $l$'s position be measured as 
%epsilon_l represents the measurement error 
%define the measured position of anchor node $l$ by
%Define the uncertain anchor positions  by
\vspace{-0.1cm}
$$
x_l=x_l^{\star}+\epsilon_l, \quad l \in \mathcal{N}_a, 
\vspace{-0.15cm}
$$
where $x_l^{\star}$ is the actual position of $l$ and 
$\epsilon_l\in\mathbb{R}^n$ represents the position information inaccuracy.
%, which is assumed as a constant.
%obeying $\|\epsilon_l\|\leq \bar{\epsilon}$, and $\bar{\epsilon}>0$ is a constant.
%Assume that  $\epsilon_l\in\mathbb{R}^n$ is a constant obeys  $\|\epsilon_l\|\leq \bar{\epsilon}$, where $\bar{\epsilon}>0$ is a constant.

%\noindent\textbf{Distance uncertainty} 
Consider the errors in the squares of inter-node distance measurements, similar to \cite{anderson2010formal}.
%Because of errors in the ranging measurements, 
%\begin{itemize}   \item
    For each  non-anchor node pair $(i,j)\in  \mathcal{E}_{ss}$, 
the square of the measured distance between them is denoted by
\vspace{-0.2cm}
$$
d_{i j}^2=d_{i j}^{\star2}+\mu_{i j},
\vspace{-0.1cm}
$$
where 
%$d_{i j}^{\star2}$ is the actual distance between them and  
$\mu_{i j}\in\mathbb{R}$ represents the measurement error.
%, which is assumed as a constant.
%obeying $\|\mu_{i j}\|\leq \bar{e}$, and $\bar{e}>0$ is a constant. 
%Assume that  $\mu_{i j}$  obeys  $\|\mu_{i j}\|\leq \bar{\mu}$, where $\bar{\mu}>0$ is a constant.
%\item 
For each anchor-non-anchor node
%anchor and non-anchor 
pair $(i,l)\in  \mathcal{E}_{as}$, the estimated distance between them is denoted by
 %with anchor uncertain 
% $$
% d_{i l}^2=d_{i l}^{\star 2}+\mu_{il}
% %\|\epsilon_l\|^2-2d_{i l}^{\star}\|\epsilon_l\| \text{cos}(\theta_{il})
% $$
% where, $\mu_{i l}\in\mathbb{R}$ small error 
% \begin{itemize}
%     \item without anchor uncertain 
%     $$
% \bar{d_{i l}}^2=d_{i l}^{\star 2}+\mu_{il}
% $$
% \item with anchor uncertain 
% $$
% d_{i l}^2=d_{i l}^{\star 2}+\|\epsilon_l\|^2-2d_{i l}^{\star}\|\epsilon_l\| \text{cos}(\theta_{il})
% $$
% \end{itemize}
%if $d_{i l}$ is only affected by anchor's uncertainty, then
$$
{
d_{i l}^2=d_{i l}^{\star 2}+{\mu}_{il},
}
\vspace{-0.1cm}
%+\mu_{il}+e_{il}}
%=d_{i l}^{\star 2}+\mu_{il}+\|\epsilon_l\|^2-2d_{i l}^{\star}\|\epsilon_l\| \text{cos}(\theta_{il})}
$$
where 
%$d_{i l}^{\star 2}$ is  the actual distance between them, 
$\mu_{i l}\in\mathbb{R}$ represents the measurement error. 
Here, $\mu_{i l}\in\mathbb{R}$ captures anchor $l$'s position  uncertainty and is assumed to be in the form
\vspace{-0.1cm}
$$
{\mu}_{il}=\|\epsilon_l\|^2-2d_{i l}^{\star}\|\epsilon_l\|\text{cos}(\theta_{il})+e_{il},
\vspace{-0.1cm}
$$
where 
$\theta_{il} \in[0,\pi]$ is the deviation angle from vector
  $x_i^{\star}-x_l^{\star}$ to vector
 %$\|x_i^{\star}-x_l^{\star}\|=d_{i l}^{\star}$ and 
% $\|x_l^{\star}-x_l\|=\|\epsilon_l\|$ 
$\epsilon_l$,
 %based on $l$'s local coordinate system,
$e_{il}\in\mathbb{R}$ is a bias term.
%, which is assumed as a constant. 
%ssume that  $e_{il}$ is a constant  obeying  $\|e_{il}\|\leq \bar{e}$, where $\bar{e}>0$ is a constant bound.
%obeying $\|\mu_{i l}\|\leq \bar{e}$, and $\bar{e}>0$ is a constant. 
It is clear that 
if anchor positions are perfectly known, then $\epsilon_l=0$ for $l\in \mathcal{N}_a$ and  hence ${\mu}_{il}=e_{il}$.

 On this basis, given the inaccurate locations $x_l$  of anchor nodes and all noisy distance measurements $d_{ij}$, $d_{il}$,
we aim to determine the locations of all non-anchor nodes and thereby formulate the noisy SNL problem:
%  our aim is to 
% { 
% %range-based SNL task in the noiseless case is to 
% determine the positions  of all non-anchor sensor nodes, $i\in \mathcal{N}_{s}$ when all anchor node positions $ x_{l} $, $l\in \mathcal{N}_{a}$ and  measurements $d_{ij}$, $d_{il}$ are given, 
% that is,  }
\vspace{-0.15cm}
\begin{align}\label{noisySNL}
	& \text{find} \quad x_{1},\dots,x_{N}\in \mathbb{R}^{n}\notag\\
	& \text{s.t.} \quad \|x_{i}-x_{j}\|^2=d_{ij}^{ 2}, \forall (i,j)\in  \mathcal{E}_{ss},\\
	& \quad\quad\;\|x_{i}-x_{l}\|^2=d_{ij}^{ 2}, \forall (i,l)\in \mathcal{E}_{as}. \notag \vspace{-0.2cm}
\end{align}
Accordingly, denote $\boldsymbol{x}^{\ddagger}=\operatorname{col}\{x_{1}^{\ddagger}, \dots, x_{N}^{\ddagger} \} \in \mathbb{R}^{nN} $ as the solution to the noisy SNL problem \eqref{noisySNL}. 

A fundamental problem in noiseless and noisy SNL problems is the existence and uniqueness of solutions  $\boldsymbol{x}^{\dagger}$ and $\boldsymbol{x}^{\ddagger}$  and their relationship with actual solution $\boldsymbol{x}^{\star}$.
%the solution $\boldsymbol{x}^{\dagger}$ or $\boldsymbol{x}^{\ddagger}$ exists and is unique. 
To this end, we formulate (\ref{noisySNL}) as a multi-player potential game to reach $\boldsymbol{x}^{\dagger}$ and $\boldsymbol{x}^{\ddagger}$ from a game-theoretic perspective.
%by reformulating .
%and  reach the solution $\boldsymbol{x}^{\ddagger}$ from
%of non-convex 
%investigate the existence condition of   %the non-convex (\ref{SNL}) as

%as the  position vector of all non-anchor nodes.

\subsection{Potential game formulation}
 { 	Here we formulate the noisy SNL  problem as an $N$-player SNL potential game $G=\{\mathcal{N}_{s}, \{\Omega_i\}_{i\in\mathcal{N}_{s}}, \{J_{i}\}_{i\in\mathcal{N}_{s}}\}$,
% as an $N$-player SNL potential game of our problem, 
 where $\mathcal{N}_{s}=\{1,\dots,N\}$ corresponds to the player set,
%, which is the set of non-anchor nodes, 
$\Omega_i$ is player $i$'s local feasible set, which is convex and compact, and $J_{i}$ is player $i$'s payoff function.  In this context, 
we map the position estimated by  each non-anchor node as each player's strategy,
i.e.,  the strategy of the player $i$ (non-anchor node) is the estimated position $x_i\in \Omega_i$.
%determined by player $i$. 
Denote $   \boldsymbol{\Omega}\triangleq\prod_{i=1}^{N}\Omega_{i} \subseteq \mathbb{R}^{nN} $, $ \boldsymbol{x}\triangleq \operatorname{col}\{x_{1}, \dots, x_{N}\} \in \boldsymbol{\Omega} $ as the position estimate strategy profile for all players,  and $ \boldsymbol{x}_{-i}\triangleq \operatorname{col}\{x_{1}, \dots,x_{i-1}, x_{i+1}, \dots, x_{N}\}\subseteq \mathbb{R}^{n(N-1)}$ as the position estimate strategy profile for all players except player $ i $. 
For $i\in\mathcal{N}_s$, the payoff function $J_{i}$ is constructed as
 \begin{equation*}\label{jI}
 J_{i}(x_{i},\boldsymbol{x}_{-i})\!=\!\sum_{j\in\mathcal{N}_{s}^{i} }(\|x_{i}-x_{j}\|^{2}\!-d_{ij}^2)^2+\!\sum_{l\in\mathcal{N}_{a}^{i} }(\|x_{i}-x_{l}\|^{2}-d_{il}^2)^2,
 %+\!\sum_{k\in\mathcal{N}_{a}^{i} }(\|x_{i}-a_{k}\|^{2}-e_{ik}^2)^2, 
  \vspace{-0.15cm}
 \end{equation*}
  where the first term in $J_{i}$  measures the   localization accuracy between  non-anchor node $i$ and its  non-anchor node neighbor $j\in\mathcal{N}_{s}^{i}$ and  the second term  measures the  localization accuracy  between  $i$ and its   anchor neighbor $l\in\mathcal{N}_{a}^{i}$.
  %and  the second term  measures the  localization accuracy  between  $i$ and its   anchor neighbor $k$. 
%这里我们考虑用最小二乘拟合优度准则来衡量节点i的定位准确性
 % Also, different from xxx, the  sensors payoffs considers the   localization information related to anchor nodes and sensor nodes.    
 % take into account the influence of anchor nodes and sensor nodes on positioning. 
% differences  denotes the sum of  differences between
%The term $\|x_{i}-x_{j}\|^{2}-d_{ij}^2$  represents the difference between  the estimated distance $\|x_{i}-x_{j}\|$  and  the  given distance $d_{ij}$, so is $\|x_{i}-a_{k}\|^{2}-e_{ik}^2$.

%To this purpose,
 The  individual objective of each non-anchor node  is to ensure its position accuracy,
 %optimal strategy $x_i$ 
 %by  considering others' strategies $\boldsymbol{x}_{-i}$, 
 i.e.,
  \vspace{-0.1cm}
 % Thus, given $\boldsymbol{x}_{-i}$,  sensor $i$ intends to solve
% adopt its location $x_{i}$ to minimize its payoff function,
\begin{equation}\label{f1}
	\min \limits_{x_{i} \in \Omega_{i}} J_{i}\left(x_{i}, \boldsymbol{x}_{-i}\right). \quad 
 \vspace{-0.15cm}
\end{equation}

Moreover, consider the following measurement of 
% the individual payoff of sensor $i$ and 
 the overall performance of sensor nodes 
% \begin{equation}\label{potential-fun}
	\begin{align}\label{potential-fun}
&
%\min _{{x}_i, i \in \mathcal{N}_s} 
{\Phi}(x_1,\dots,x_{{N}})\\
&=  \sum_{(i, j) \in \mathcal{E}_{ss}} \!\!\!\!{(\left\|{x}_i-{x}_j\right\|^2-d_{i j}^2)^2}
%+  \sum_{l \in \mathcal{N}_a}\!\!\left\|{x}_l-x_l\right\|^2
+\sum_{(i, l) \in \mathcal{E}_{as}} \!\!\!\!{(\left\|{x}_i-{x}_l\right\|^2-d_{i l}^2)^2} 
%+  \sum_{l \in \mathcal{N}_a}\!\!\left\|{x}_l-x_l\right\|^2
\notag\\
&=  \!\sum_{(i, j) \in \mathcal{E}_{ss}} \!\!\!\!{(\left\|{x}_i\!-\!{x}_j\right\|^2\!-\!(d_{i j}^{\star 2}\!+\!\mu_{i j}))^2}\!\notag\\
%+\!  \sum_{(i, l) \in \mathcal{E}_{as}} \!\!\!\!{(\left\|{x}_i\!-\!(x_l^{\star}\!+\!\epsilon_l)\right\|^2\!-\!(d_{i l}^{\star 2}\!+\!\mu_{il}}\\
& \quad+\!\!\!\!\!\!\!\sum_{(i, l) \in \mathcal{E}_{as}} \!\!\!\!\!\!(\left\|{x}_i\!-\!(x_l^{\star}\!+\!\epsilon_l)\right\|^2\!\!\!-\!(d_{i l}^{\star 2}\!+\!\!\|\epsilon_l\|^2\!\!-\!\!2d_{i l}^{\star}\|\epsilon_l\|\text{cos}(\theta_{il})\!\!+\!\!e_{il}))^2. \notag
% &=\sum_{(i, j) \in \mathcal{E}\backslash\mathcal{E}_{aa}} \!\!\!\!{(\left\|{x}_i-{x}_j\right\|^2-(d_{i j}^{\star 2}+\mu_{i j}))^2}+\|\Lambda(\boldsymbol{x}-(\boldsymbol{x}^{\star}+\boldsymbol{\epsilon}))\|^2
 \vspace{-0.15cm}
\end{align}
%\end{equation}
{Here,  $J_{i}$ denotes the   localization accuracy of  node $i$, which depends on the  strategies of $i$'s neighbors, while  $	\Phi$ denotes the   localization accuracy of the entire network $\mathcal{G}$.

%\subsection{Potential game model to SNL problem}\label{for-potential}
%\subsection{Preliminary of potential game}\label{pre-potential}
%直接放在开头写建模成势博弈

%\subsection{Potential game formulation}
From a game-theoretic perspective, each sensing agent is considered as a selfish entity who simply tries to minimize its own payoff function. For this $N$-player game to provide a solution to localize the whole sensor network,
%In the  SNL problem,
each non-anchor node needs to consider the location accuracy of the whole sensor network while ensuring its own positioning accuracy through the given information.
%To achieve this alignment,
%In this section,
%To  better describe  sensor nodes' individual preference and network interaction,
% In other words, each non-anchor node needs to guarantee consistency between its  individual objective and   collective objective.
% %, which will be defined explicitly in the sequel. 
% %To this end, by regarding the individual payoff  $J_{i}$ as a marginal contribution to the whole network's collective objective \cite{marden2009cooperative,jia2013distributed},
% we consider the following measurement of 
% % the individual payoff of sensor $i$ and 
%  the overall performance of sensor nodes %Here, we employ  the Least-Squares goodness of fit criterion  to measure the localization accuracy of sensor node $i$. 
% \begin{equation}\label{potential-fun}
% 	\Phi\!\left(x_{1},\!\dots,\! x_{N}\right)\!=\sum_{(i,j)\in\mathcal{E} }(\|x_{i}-x_{j}\|^{2}\!-d_{ij}^2)^2.
%  %\!+\!\!\!\!\!\!\sum_{(i,k)\in\mathcal{E}_{as} }\!\!\!\!\!(\|x_{i}-a_{k}\|^{2}-e_{ik}^2)^2.
% \end{equation}
% {Here,  $J_{i}$ denotes the   localization accuracy of non-anchor node $i$, which depends on the  strategies of $i$'s neighbors, while  $	P$ denotes the   localization accuracy of the entire network $\mathcal{G}$. 
On this basis,   we  formulate $G=\{\mathcal{N}_{s}, \{\Omega_i\}_{i\in\mathcal{N}_{s}}, \{J_{i}\}_{i\in\mathcal{N}_{s}}\}$ as a potential	game,  where  $\Phi$ in (\ref{potential-fun}) satisfies the concept of a potential function \cite{monderer1996potential}, i.e.,
 \vspace{-0.1cm}
\begin{equation}\label{pp1}
	\Phi(x_{i}^{\prime}, \boldsymbol{x}_{-i})-\Phi\left(x_{i}, \boldsymbol{x}_{-i}\right)=J_{i}(x_{i}^{\prime}, \boldsymbol{x}_{-i})-J_{i}\left(x_{i}, \boldsymbol{x}_{-i}\right),
  \vspace{-0.15cm}
	\end{equation}
for  every $i\in\mathcal{N}_s$,  $\boldsymbol{x}\in \boldsymbol{\Omega}$,  and unilateral deviation $x_{i}^{\prime}\in \Omega_i$.
% \begin{mydef}[potential game]\label{d2}
% 	%From the network’s perspective, it is desirable to obtain the optimum location selection  which minimizes the network cost function given by
% 	A game $G=\{\mathcal{N}_{s}, \{\Omega_i\}_{i\in\mathcal{N}_{s}}, \{J_{i}\}_{i\in\mathcal{N}_{s}}\}$ is a potential
% 	game if there exists a potential function $\Phi$ such that, for $i\in\mathcal{N}_s$,
%   \vspace{-0.1cm}
% 	\begin{equation}\label{pp1}
% 	\Phi(x_{i}^{\prime}, \boldsymbol{x}_{-i})-\Phi\left(x_{i}, \boldsymbol{x}_{-i}\right)=J_{i}(x_{i}^{\prime}, \boldsymbol{x}_{-i})-J_{i}\left(x_{i}, \boldsymbol{x}_{-i}\right),
%   \vspace{-0.15cm}
% 	\end{equation}
% 	for every  $\boldsymbol{x}\in \boldsymbol{\Omega}$,  and unilateral deviation $x_{i}^{\prime}\in \Omega_i$. 
% \end{mydef}
% It follows from Definition   \ref{d2} that 
% %the change in each player’s payoff is influenced by the change of a global unified  potential function.
% 
This indicates that any unilateral deviation from a strategy profile always results in the same change in both individual payoffs and a unified potential function.  In other words, 
%In order for thisN-person game to provide a solution to the cooperative sensing problem in (P), 
the  individual goal  $J_{i}$  is aligned with the global objective $	\Phi$.
Moreover, to attain an optimal value for $J_{i}\left(x_{i}, \boldsymbol{x}_{-i}\right)$,   players need to engage in negotiations and alter their optimal strategies. The best-known concept that describes an acceptable result achieved by all players is  NE \cite{nash1951non}.
%,  whose definition is formulated below. 
\begin{mydef}[NE]\label{d1}
	A profile $ \boldsymbol{x}^{\Diamond}=\operatorname{col}\{x_{1}^{\Diamond}, \dots, x_{N}^{\Diamond} \} \in \boldsymbol{\Omega} \subseteq \mathbb{R}^{nN}$ is said to be an NE of game (\ref{f1}) if for any $i \in \mathcal{N}_s  $,   
 %for any $ x_{i}\in \Omega_{i}$  we have
	\begin{equation}\label{ne}
		J_{i}\left(x_{i}^{\Diamond}, \boldsymbol{x}_{-i}^{\Diamond}\right) \leq J_{i}\left(x_{i}, \boldsymbol{x}_{-i}^{\Diamond}\right),\;\forall x_{i}\in \Omega_{i}.
	\end{equation}
\end{mydef}
%The NE above characterizes a strategy profile $\boldsymbol{x}^{\star}$ that  each non-anchor adopts its optimal position. Given others' strategies, no non-anchor node can benefit from changing its position unilaterally. 
%It follows from Definition \ref{d2}
%the potential function can ensure the optimization of local sensor location and the improvement of the whole network performance. 
% Note that an NE of a potential game ensures not only  that each non-anchor node can adopt its optimal location strategy from the individual perspective, but also that  the sensor network as a whole can achieve a precise localization  from the global perspective.
% , i.e., 
%  %. That is to say, the global NE is specified as 
%  the strategy profile  $\boldsymbol{x}^{\star}$ satisfies $\|x_{i}^{\star}-x_{j}^{\star}\|^2-d_{i j}^2=0$ for any  $(i,j)\in \mathcal{E}$.
Here, we call  NE as \textit{global} NE due to the non-convex SNL formulation in this paper. This is different from the concept of \textit{local} NE \cite{pang2011nonconvex,heusel2017gans}, which  only satisfies condition \eqref{ne} within a small neighborhood of $x_{i}^{\Diamond}$ for $i \in \mathcal{N}_s$, rather than  the whole $\Omega_{i}$. 
% We also  consider another mild but
% well-known concept 
% to help characterize the solutions to (\ref{f1}).
\begin{mydef}[local NE]\label{ce234}	
	A strategy profile $ \boldsymbol{x}^{\Diamond} $ is said to be a local NE  of  (\ref{f1}) if there exists a constant $r > 0$ such that for any $i \in \mathcal{N}_s  $,  
	%the point $x_i^{\star}$ satisfies 
	% A strategy profile $ \boldsymbol{x}^{\star} $ is said to be a \textbf{local Nash equilibrium} of game (\ref{f1}) if, for all $ i\in \mathcal{I}$, $x_{i}^{\star}\in \Omega_i$, $(x_{i}^{\star},\boldsymbol{x}_{-i}^{\star})$ satisfying 
	% 		\begin{equation}\label{ee1}
	% 	J_{i}(x_{i}^{\star}, \boldsymbol{x}_{-i}^{\star})\leq J_{i}(x_{i}, \boldsymbol{x}_{-i}^{\star}), \quad \forall x_{i}\in \Omega_i \cap   B_{\delta}(x_{i}^{\star})
	% 		\end{equation}
	\begin{equation}\label{e234}
	\begin{aligned}
	J_{i}\left(x_{i}^{\Diamond}, \boldsymbol{x}_{-i}^{\Diamond}\right) \leq J_{i}\left(x_{i}, \boldsymbol{x}_{-i}^{\Diamond}\right),\;\forall x_i \in \Omega_i \cap   \mathcal{B}_{r}(x_{i}^{\Diamond}).
	%	&\mathbf{0}_{m} \in-g_{i}\left(x_{i}^{\star}\right)+\mathcal{N}_{\mathbb{R}_{\geq 0}}^{m}\left(\lambda_{i}^{\star}\right)
	\end{aligned}
	\end{equation}
 where $\mathcal{B}_{r}(x_{i}^{\Diamond})\triangleq \{y\in \mathbb{R}^n : \|y-x_{i}^{\Diamond}\| < r\}$ is an open  Euclidean ball with radius $r$ and center $x_{i}^{\Diamond}$. 
 %{where $\mathcal{N}_{\Omega_{i}}(x_{i}^{\star})=\{e\in\mathbb R^n:e^T(x-x_{i}^{\star})\leq 0,\forall x\in\Omega_{i}\}$ is the normal cone at point $x_{i}^{\star}$ on set $\Omega_{i}$.}
	%	with $\boldsymbol{x}_{-i}=\boldsymbol{x}_{-i}^{\star}$.
\end{mydef}
% It is not difficult to reveal that in non-convex games, if $ \boldsymbol{x}^{\star} $ is a global NE, then it must be a NE stationary point, but not vice versa.
% % Here we call  NE as global NE in this paper due to the non-convex formulation \eqref{f1}, which is different from  the concept of \textit{local} NE \cite{nouiehed2019solving,heusel2017gans} that only satisfying \eqref{ne} on a small neighbor of $x_{i}^{\star}$ for $i \in \mathcal{N}_s$, rather than the whole $\Omega_{i}$.

In order to reveal the relations between  NE $\boldsymbol{x}^{\Diamond}$, noiseless node positions $\boldsymbol{x}^{\dagger}$ and noisy node positions $\boldsymbol{x}^{\ddagger}$, we need to utilize graph rigidity theory \cite{anderson2008rigid}, particularly the notions and tools summarized in Appendix \ref{lea1}.
%investigate the existence and uniqueness of  NE 
%in the noisy SNL problem, 
With these notions, we make the following generic and feasible assumption. 

\begin{assumption}
%	$\ $
	%\begin{itemize}
		%	\item The range measurements $d_{ij}$ and $e_{ik}$ are noisy-free and 		symmetric, and all anchors' positions are accurate.
			% For $ i \in \mathcal{I} $, $ \Omega_{i} $ is compact and convex.	
	%	\item 
 The sensor topology graph	 $ \mathcal{G} $ is undirected and generically globally rigid.
	%\end{itemize}
\end{assumption}
%The  assumptions of convexity and compactness about local feasible are quite common and {have been} widely used in the related wireless sensor network literature \cite{jia2013distributed,ke2017distributed}. 
% The undirected graph topology is a common assumption in many graph-based approaches \cite{chen2021distributed,jing2021angle}. 
% The connectivity of $\mathcal{G} $ can also be induced by some disk graph \cite{wan2019sensor},  which ensures the validity of the information transmission between nodes.
%我们的东西也是可以从disk graph诱导，我们提供connected,来保证节点之间的信息可以有效传递，
%我们和disk graph可以互相融合
%The assumption about noisy-free and symmetric measurements and  position accuracy of anchors  is quite common  in the related wireless sensor network literature \cite{jia2013distributed,wan2019sensor}. 
%the undirected graph is usually a common assumption in many graph-based problems \cite{chen2021distributed,jing2021angle}. 
{The generic global rigidity of $\mathcal{G} $ has been widely employed in SNL problems without measurement noises to make the geometric realization of the graph invariant, which indicates unique localization of the sensor network \cite{eren2004rigidity,anderson2008rigid}. 
%That is to say, we can establish $\boldsymbol{x}^{\star}=\boldsymbol{x}^{\dagger}$ in the noiseless case \cite{eren2004rigidity}. Moverover, it follows from \cite{anderson2010formal,xu2024global} that if all $\mu_{ij}$, $e_{il}$ and $\epsilon_l$ are zero, Assumption 1 guarantees that the  NE  $\boldsymbol{x}^{\Diamond}$ of  potential game ${{G}}$ is unique and equal to the solution $\boldsymbol{x}^{\star}$.
}
 % and has raised extensive discussions in existing works \cite{wan2019sensor,xxx}
%Besides, there have been extensive discussions on graph rigidity in existing works \cite{wan2019sensor,cao2021bearing}, but it is not the primary focus of our paper. }

%Moreover, it follows from \cite{anderson2010formal,xu2024global} that if all $\mu_{ij}$, $e_{il}$ and $\epsilon_l$ are zero the  NE  $\boldsymbol{x}^{\Diamond}$ of  potential game ${{G}}$ is unique and equal to the solution $\boldsymbol{x}^{\star}$  with . 
%If all $\mu_{ij}$, $e_{il}$ and $\epsilon_l$ are zero, noisy SNL problem \eqref{noisySNL} is converted into noiseless problem \eqref{SNL}. It follows from \cite{anderson2010formal,xu2024global} that  the  NE  $\boldsymbol{x}^{\Diamond}$ of  potential game ${{G}}$ is unique and equal to the solution $\boldsymbol{x}^{\star}$  of \eqref{SNL}.
\subsection{Existence, uniqueness, and errors of the solution}
As for the noiseless case, the global rigidity of the sensor network graph in Assumption 1 guarantees that (\ref{SNL}) has a unique solution $\boldsymbol{x}^{\dagger}$  equal to the  actual positions $\boldsymbol{x}^{\star}$ \cite{eren2004rigidity}. We uniformly use $\boldsymbol{x}^{\star}$ to represent $\boldsymbol{x}^{\dagger}$ hereafter. 
Moverover,  we establish in the next section that if all $\mu_{ij}=0$ and $e_{il}=0$, the  NE  $\boldsymbol{x}^{\Diamond}$ of  potential game ${{G}}$ is unique and equal to $\boldsymbol{x}^{\star}$.
% with all $\mu_{ij}=0$ and $e_{il}=0$ and $\epsilon_l$ are zero
% it follows from \cite{anderson2010formal,xu2024global} that if all $\mu_{ij}$, $e_{il}$ and $\epsilon_l$ are zero, Assumption 1 guarantees that the  NE  $\boldsymbol{x}^{\Diamond}$ of  potential game ${{G}}$ is unique and equal to the solution $\boldsymbol{x}^{\star}$.
%(\ref{SNL}) has a unique solution $\boldsymbol{x}^{\dagger}$ and is equal to $\boldsymbol{x}^{\star}$ \cite{eren2004rigidity}. Moreover,  it follows from \cite{anderson2010formal,xu2024global} that  the  NE  $\boldsymbol{x}^{\Diamond}$ of  potential game ${{G}}$ is unique and corresponding to the a $\boldsymbol{x}^{\star}$  of \eqref{SNL}.
%Assumption 1 guarantees that (\ref{SNL}) has a unique solution $\boldsymbol{x}^{\dagger}$ and is equal to $\boldsymbol{x}^{\star}$ \cite{eren2004rigidity}.
%the global rigidity of the sensor network graph guarantees the existence and uniqueness of  NE $\boldsymbol{x}^{\Diamond}$ in the SNL problem, which is  equal to the actual sensor positions, i.e,  $\boldsymbol{x}^{\Diamond}=\boldsymbol{x}^{\dagger}=\boldsymbol{x}^{\ddagger}=\boldsymbol{x}^{\star}$.
%there exists a unique NE of the SNL problem, equal to the actual sensor positions $\boldsymbol{x}^{\star}$. 

However, in the noisy case, obtaining the relationship among $\boldsymbol{x}^{\Diamond}$, $\boldsymbol{x}^{\ddagger}$, and $\boldsymbol{x}^{\star}$ is not as straightforward. The presence of noise may perturb the rigidity of the graph, resulting in flip ambiguities or even no solution to the problem.  Fig. \ref{fig123s}(a) considers a case with inaccuracies in anchor node position information.  
%Fig. \ref{fig123s}(a) depicts a sensor network with two non-anchor nodes and three anchor nodes in two configurations.
The red pluses denote the true anchor node locations, while the blue stars represent the true non-anchor node locations. The grey lines indicate connections in this configuration. On the other hand, 
the red circles denote the noisy anchor node location information, while 
the green circles denote the computed noisy non-anchor node locations.
% the noisy anchor nodes are shown by red circles and the computed noisy non-anchor nodes are shown by green circles. 
The dashed lines indicate the noisy inter-node distance measurements. It can be seen that both anchor positions and inter-node distances among nodes in the actual configuration are close to those in the noisy one. However, the calculated or estimated 
 positions of the non-anchor nodes become flipped compared to the actual positions.
%If there exists a bound such that when these measurement noises do not exceed it, then such a flipping case would not occur. 
Furthermore,  consider another case with errors in both anchor node position information and inter-node distance measurements. As shown in Fig. \ref{fig123s}(b),  even small errors may lead to the final localization results deviating largely from the actual positions.  
%measures measurement error, 
%such that

% . The true
% locations of non-anchor nodes are shown by blue diamonds.

% , one corresponding to nodes 1, 2, 3,
% and 4 and the other corresponding to nodes . Suppose that nodes 2, 3,
% and 4 are anchors, and regard 1 as being in a true position. It is quite apparent from
% the figure that the distances d12, d13, d14 are close to the distances d12, d13, d14.
	\begin{figure}[ht]
			\hspace{-0.7cm}
			\centering	
			\subfigure[ ]{
				\begin{minipage}[t]{0.48\linewidth}
					\centering
					\includegraphics[width=4.6cm]{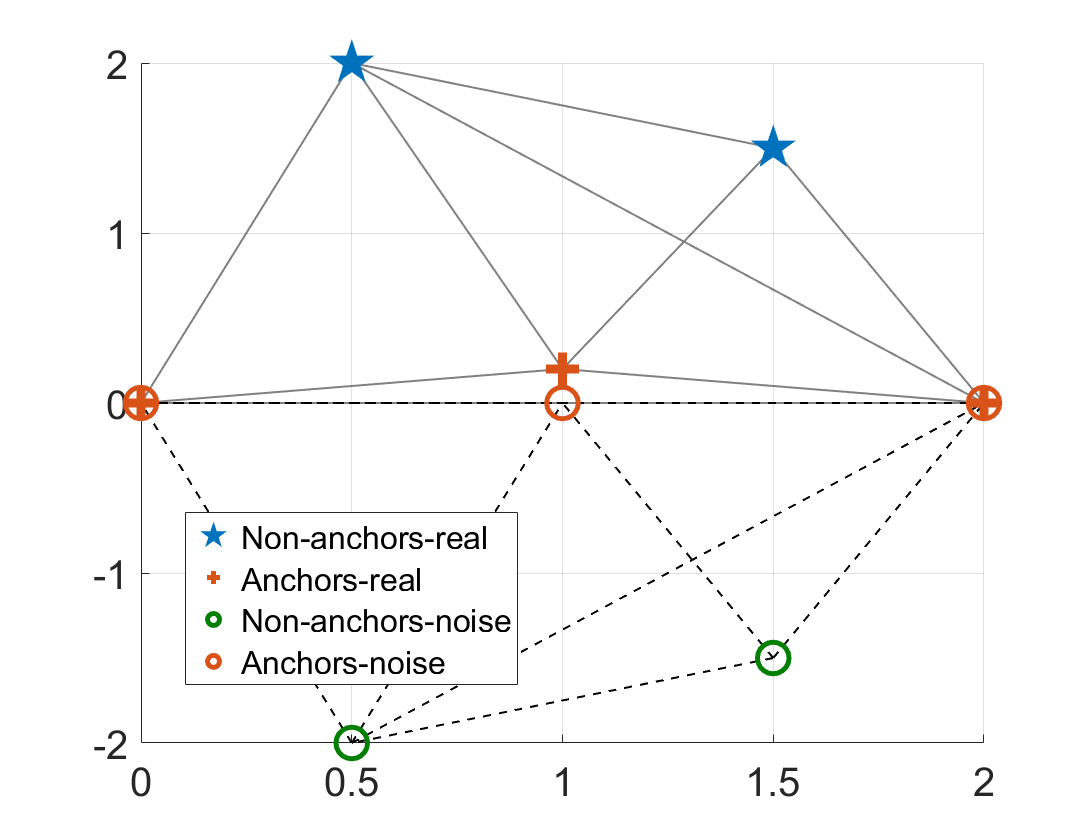}
					%\caption{fig1}
				\end{minipage}%
			}%
			\hspace{0.1cm}
			\subfigure[ ]{
				\begin{minipage}[t]{0.48\linewidth}
					\centering
					\includegraphics[width=4.6cm]{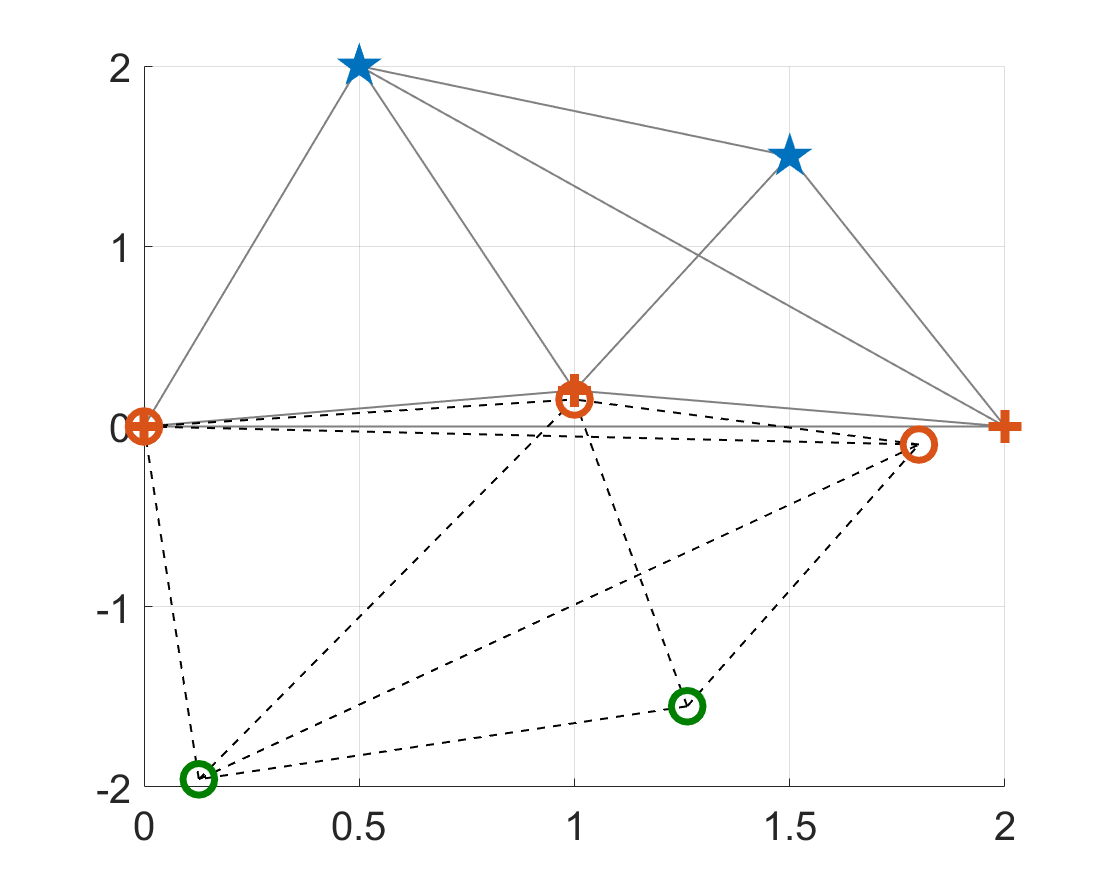}
					%\caption{fig2}
				\end{minipage}%
			}%
			\centering
   \vspace{-0.15cm}
			\caption{A  sensor network  with two non-anchor nodes and three anchor nodes in two configurations.} 
   \vspace{-0.3cm}
			\label{fig123s}
		\end{figure}
  
%  \begin{figure}[ht]
% 	\centering	
% 	\includegraphics[scale=0.2]{n5.png}\\
%  \vspace{-0.2cm}
% 	\caption{A sensor network  with two non-anchor nodes and three anchor nodes in two configurations.}
% 	\label{fig123s}
% 	\vspace{-0.3cm}
% \end{figure}

%  \begin{figure}[ht]
% 	\centering	
% 	\includegraphics[scale=0.2]{noise4.png}\\
%  \vspace{-0.2cm}
% 	\caption{A sensor network  with two non-anchor nodes and three anchor nodes in two configurations.}
% 	\label{fig1234s}
% 	\vspace{-0.3cm}
% \end{figure}
%the disturbance of noise
%in the presence of 
%In this paper, we aim to investigate the existence and uniqueness of  NE in the noisy SNL problem, and reveal its relation with the true sensor positions.  Unlike the noiseless case where the existence and uniqueness of NE can be obtained based on the global rigidity of the sensor network graph \cite{anderson2010formal,xu2024global}, the disturbance of noise may break the rigidity of the graph, resulting in no solution to the problem. 
%Therefore, we consider under what level of error solutions exist and do not deviate too much from the true solution."
%我们考虑在哪种误差程度下有解。
Instead, consider a threshold for measurement noises such that if the errors in anchor node position and inter-node distance measurements do not exceed it, the extreme cases mentioned above would not occur and the node positions would remain close to the correct values. Therefore,
%不同于无噪声的情况下NE和存在性唯一性可以直接由传感器网络图的全局刚性得到， 噪声的存在可能会破坏图的刚性，从而导致问题无解。
the questions of interest to us are  as follows:
%改成两个点
%expect to solve the following problems
%\vspace{-0.15cm}
%\begin{problem}
%$~$
 \begin{enumerate}[(i)]
 %	\item  
 	%In what conditions, the stationary point coincide with the global NE?
 	%	\item if such an optimal solution exists, whether it is also the precise localization of the whole network,
 %	\item whether such an optimal solution  is also the precise localization of the whole network,
 	\item Under what conditions on noise levels does there exist a unique  NE  $\boldsymbol{x}^{\Diamond}$ and when does it correspond to the solution of the SNL problem \eqref{noisySNL}?

%  is equal to the solution $\boldsymbol{x}^{\ddagger}$ of noisy SNL problem \eqref{noisySNL}?
  %Under what conditions of noise levels does there exist an  NE  $\boldsymbol{x}^{\Diamond}$?
 
%  Is it guaranteed that, at least for suitably small noise, there will still be an NE to the SNL problem and that its deviation from the true sensor positions will go continuously to zero as the noise goes to zero?
 	% \item 
  % %Under what conditions of noise levels,
  % Under what conditions of noise levels will the NE $\boldsymbol{x}^{\Diamond}$  be the unique solution  $\boldsymbol{x}^{\ddagger}$ to noisy SNL problem \eqref{noisySNL}?
  \item How do the errors between actual positions $\boldsymbol{x}^{\star}$ and NE  $\boldsymbol{x}^{\Diamond}$ depend on the measurement noise amplitudes?
  %the solution $\boldsymbol{x}^{\ddagger}$ of the noisy SNL problem 
  %\item How to design an efficient algorithm to seek the  NE $\boldsymbol{x}^{\Diamond}$?
 \end{enumerate}
%\end{problem}
  
\section{Existence and Uniqueness of NE}

%xxxxxx 

%Next, we show that global NE $\boldsymbol{x}^{\star}$ is unique and  represents the actual position profile of all non-anchor nodes in the noiseless case, which is equal to the global solution of the SNL.

In this section, we investigate the existence and uniqueness of global NE $\boldsymbol{x}^{\Diamond}$. 
%We first consider the noiseless case.
For the noiseless case, the following theorem reveals that the  global NE $\boldsymbol{x}^{\Diamond}$ exists and is unique, 
corresponding to the actual non-anchor nodes' positions $\boldsymbol{x}^{\star}$.
The proof is given in Appendix \ref{rea2}.

{
\begin{mythm}\label{l11} 
% These equations have the solution $x_i =x_i^{\star}$, $i\in\mathcal{N}_s$ and $x_l =x_l^{\star}$, $l\in\mathcal{N}_a$ when the $\mu_{ij}=0$ for
% $(i,j) \in\mathcal{E}\backslash \mathcal{E}_{aa}$ and $\epsilon_l=0$ for $l\in\mathcal{N}_a$ .
Under Assumption 1, there exists a unique NE $\boldsymbol{x}^{\Diamond}$ of the game $G$, which satisfies  $x_i^{\Diamond} =x_i^{\star}$ for $i\in\mathcal{N}_s$ and $x_l =x_l^{\star}$ for $l\in\mathcal{N}_a$ if $\mu_{ij}=0$ for
$(i,j) \in\mathcal{E}_{ss} $, ${\mu}_{il}=0$ for
$(i,j) \in\mathcal{E}_{as} $ and $\epsilon_l=0$ for $l\in\mathcal{N}_a$.
\end{mythm}
 %the existence and uniqueness of NE in the noiseless SNL problem are 
  % Theorem \ref{l11} is guaranteed by the global rigidity of the sensor network graph. However, the global rigidity may be disrupted 
  On the other hand, when the measurement noises arise, we first investigate under what conditions of noise levels a unique local NE exists and is near the true sensor node locations.
%the existence of a unique local NE for 
Denote $\boldsymbol{\mu}=\text{col}\{\mu_{ij}\}_{(i,j)\in \mathcal{E}_{ss}}$, %$\bar{\boldsymbol{e}}_2=\text{col}\{\bar{\mu}_{il}\}_{(i,l)\in \mathcal{E}_{as}}$,
$\boldsymbol{e}=\text{col}\{e_{il}\}_{(i,l)\in \mathcal{E}_{as}}$, $\boldsymbol{\epsilon}=\text{col}\{\epsilon_l\}_{l\in\mathcal{N}_a}$, $\boldsymbol{\mu}_0=\text{col}\{\mu_{ij}=0\}_{(i,j)\in \mathcal{E}_{ss}}$, 
$\boldsymbol{e}_0=\text{col}\{e_{il}=0\}_{(i,l)\in \mathcal{E}_{as}}$, $\boldsymbol{\epsilon}_0=\text{col}\{\epsilon_l=0\}_{l\in\mathcal{N}_a}$.

%The following lemma establishes that there exists a local NE that is unique.
The following lemma is established for the SNL problem with measurement noise. The proof is given in Appendix \ref{rea}.
%, provided that  the errors are small enough.
\begin{lemma}\label{ee1}
 $~$
%Under Assumption 1, consider the set of equations$\nabla \Phi=0$, which are necessarily satisfied at a minimum of \eqref{robust1}. Then the following are true:
\begin{enumerate}[i)]
%     \item These equations have the solution $x_i =x_i^{\star}$, $i\in\mathcal{N}_s$ and $x_l =x_l^{\star}$, $l\in\mathcal{N}_a$ when the $\mu_{ij}=0$ for
% $(i,j) \in\mathcal{E}\backslash \mathcal{E}_{aa}$ and $\epsilon_l=0$ for $l\in\mathcal{N}_a$ .
\item {
%Given fixed $\theta_{il}\in [0,\pi]$ for $l\in\mathcal{N}_a$,  $\boldsymbol{\mu}$, $\boldsymbol{e}$, and $\boldsymbol{\epsilon}$.  
Under Assumption 1, there exists a small positive $\delta$ and positive constants $a$,  $b$ and $c$ depending on $\delta$ such that  for any constant uncertainty vectors $\boldsymbol{\mu}$, $\boldsymbol{e}$, and $\boldsymbol{\epsilon}$ satisfying $a\|\boldsymbol{\mu}\|^2+b\|\boldsymbol{e}\|^2+c\|\boldsymbol{\epsilon}\|^2< \delta$,  
%suitably small positive $\delta$, if $a\|\boldsymbol{\mu}_1\|^2+b\|\boldsymbol{\epsilon}\|^2+c\|\boldsymbol{e}_2\|^2< \delta$, 
% Given fixed
% %$\theta_{il}\in [0,\pi]$ for $l\in\mathcal{N}_a$, 
% $\boldsymbol{\mu}_1$, $\bar{\boldsymbol{\mu}}_2$, and $\boldsymbol{\epsilon}$.  There exist    positive constants $a$,  $b$, $c$ and $d$  such that for any
% suitably small positive $\delta$, if $a\|\boldsymbol{\mu}_1\|^2+b\|\bar{\boldsymbol{\mu}}_2\|^2+c\|\boldsymbol{\epsilon}\|^2+d\|\boldsymbol{\epsilon}\|^4< \delta$, 
% There exist  suitably small positive $\delta$,
% %$\delta_1$, $\delta_2$, $\delta_3$ 
%  given fixed $\theta_{il}\in [0,\pi]$ for $l\in\mathcal{N}_a$
% and  positive constants $a$,  $b$ and $c$ depending on $\delta$
% %, $\delta_2$ and $\delta_3$
% such that for any fixed $\boldsymbol{\mu}_1$, $\boldsymbol{\mu}_2$, and $\boldsymbol{\epsilon}$ satisfying $a\|\boldsymbol{\mu}_1\|^2+b\|\boldsymbol{\mu}_2\|^2+c\|\boldsymbol{\epsilon}\|^2< \delta$, 
%lying in the ball of radius $\delta_1$, and any fixed $\boldsymbol{\mu}_2$ lying in the ball of radius $\delta_2$, any fixed $\boldsymbol{\epsilon}$ lying in the ball of radius $\delta_3$ around the origin, 
 there is a unique solution  $ \hat{\boldsymbol{x}}=\operatorname{col}\{\hat{x}_{1}, \dots, \hat{x}_{N} \} \in \boldsymbol{\Omega} \subseteq \mathbb{R}^{nN}$ of the equations $\nabla {\Phi}=0$
 %of the equations $\nabla \Phi=0$
satisfying  
$\|\hat{\boldsymbol{x}}-\boldsymbol{x}^{\star}\|^2< \delta$.
%${\Phi}(\hat{\boldsymbol{x}})-{\Phi}(\boldsymbol{x}^{\star})< \delta$.
% $\|\hat{\boldsymbol{x}}-\boldsymbol{x}^{\star}\|\rightarrow 0$ as $\|\boldsymbol{\mu}_1\|\rightarrow 0$, $\|\boldsymbol{\mu}_2\|\rightarrow 0$ and $\|\boldsymbol{\epsilon}\|\rightarrow 0$.
}
% the constraint $\|\hat{\boldsymbol{x}}-\boldsymbol{x}^{\star}\|\leq c\|\boldsymbol{\mu}\|+d\|\boldsymbol{\epsilon}\| $
% or $\|\hat{\boldsymbol{x}}-\boldsymbol{x}^{\star}\|\leq f(\|\boldsymbol{\mu}\|,\|\boldsymbol{\epsilon}\|) $, where
\item  Under Assumption 1, for constant uncertainty vectors $\boldsymbol{\mu}$, $\boldsymbol{e}$ and $\boldsymbol{\epsilon}$, the solution $ \hat{\boldsymbol{x}}$ is a local NE (with respect to $\boldsymbol{x}$) of
${\Phi}(\boldsymbol{x},\boldsymbol{\mu},{\boldsymbol{e}},\boldsymbol{\epsilon})$.
\end{enumerate}   
\end{lemma}

Lemma \ref{ee1} indicates that small magnitudes of measurement noises and uncertainty vectors will not disrupt the global rigidity of the graph, i.e., the positive definiteness of the Hessian matrix $\nabla^2 {\Phi}$. Accordingly, the deviation of a local NE  from the actual sensor positions will be bounded and continuously converge to zero as the errors in anchor node position information and inter-node distance measurements approach zero.

Next, 
let $\mathcal{B}_{\delta_1}(\boldsymbol{x}^{\star})$ denote the ball around $\boldsymbol{x}^{\star}$ defined by $\|\boldsymbol{x}-\boldsymbol{x}^{\star}\|^2< \delta_1$, where $\delta_1$ is a  positive constant.
 Let $\mathcal{B}^c_{\delta_1}(\boldsymbol{x}^{\star})$ denote the complementary set $\|\boldsymbol{x}-\boldsymbol{x}^{\star}\|^2\geq \delta_1$. Observe that for all $\boldsymbol{\mu}$, $\boldsymbol{e}$ and $\boldsymbol{\epsilon}$
with $a\|\boldsymbol{\mu}\|^2+b\|\boldsymbol{e}\|^2+c\|\boldsymbol{\epsilon}\|^2< \delta_1$, Lemma \ref{ee1} guarantees that $\hat{\boldsymbol{x}} \in \mathcal{B}_{\delta_1}(\boldsymbol{x}^{\star})$. Let ${\Phi}_1$
be defined by
\vspace{-0.1cm}
\begin{align}\label{p11}
{\Phi}_1=&\min _{{x}_i, i \in \mathcal{N}_s, \boldsymbol{x} \in \mathcal{B}^c_{\delta_1}\!(\boldsymbol{x}^{\star})} 
%{\Phi}_1(x_1,\dots,x_{{N}})\\&= 
\sum_{(i, j) \in \mathcal{E}_{ss}} \!\!\!\!{(\left\|{x}_i-{x}_j\right\|^2-d_{i j}^{\star2})^2}\\
%+  \sum_{l \in \mathcal{N}_a}\!\!\left\|{x}_l-x_l\right\|^2
&\quad\quad\quad\quad\quad\quad\quad+\!\!\!\!\sum_{(i, l) \in \mathcal{E}_{as}} \!\!\!\!{(\left\|{x}_i-{x}_l^{\star}\right\|^2-d_{i l}^{\star2})^2}. \notag
%+  \sum_{l \in \mathcal{N}_a}\!\!\left\|{x}_l-x_l\right\|^2
% &=\sum_{(i, j) \in \mathcal{E}\backslash\mathcal{E}_{aa}} \!\!\!\!{(\left\|{x}_i-{x}_j\right\|^2-(d_{i j}^{\star 2}+\mu_{i j}))^2}+\|\Lambda(\boldsymbol{x}-(\boldsymbol{x}^{\star}+\boldsymbol{\epsilon}))\|^2
\vspace{-0.35cm}
\end{align}

On the other hand,  consider also a collection of minimization problems, parameterized by a nonnegative constant $\delta_2$, with variables $\boldsymbol{x}$, $\boldsymbol{\mu}$, $\boldsymbol{e}$ and $\boldsymbol{\epsilon}$:
\begin{align}\label{phi2}
\Phi_2=\!\!\!\!\!\!\!\!&\min _{{x}_i, i \in \mathcal{N}_s, \boldsymbol{x} \in \mathcal{B}^c_{\delta_1}\!(\boldsymbol{x}^{\star}), a\|\boldsymbol{\mu}\|^2+b\|\boldsymbol{e}\|^2+c\|\boldsymbol{\epsilon}\|^2\leq \delta_2} %{\Phi}_2(x_1,\dots,x_{{N}})\\&
\sum_{(i, j) \in \mathcal{E}_{ss}} \!\!\!\!(\left\|{x}_i-{x}_j\right\|^2\notag\\
%+  \sum_{l \in \mathcal{N}_a}\!\!\left\|{x}_l-x_l\right\|^2
&\quad-d_{i j}^{\star2}-\mu_{ij})^2+\!\!\!\!\!\!\!\sum_{(i, l) \in \mathcal{E}_{as}} \!\!\!\!\!\!(\left\|{x}_i\!-\!(x_l^{\star}\!+\!\epsilon_l)\right\|^2\!\!\!-\!(d_{i l}^{\star 2}\!+\!\|\epsilon_l\|^2 \notag\\
&\quad-\!\!2d_{i l}^{\star}\|\epsilon_l\|\text{cos}(\theta_{il})\!\!+\!\!e_{il}))^2. 
\vspace{-0.35cm}
%+\sum_{(i, l) \in \mathcal{E}_{as}} \!\!\!\!{(\left\|{x}_i-{x}_l^{\star}-\epsilon_l\right\|^2-d_{i l}^{\star2}-\mu_{il}-e_{il})^2} \notag
%+  \sum_{l \in \mathcal{N}_a}\!\!\left\|{x}_l-x_l\right\|^2
% &=\sum_{(i, j) \in \mathcal{E}\backslash\mathcal{E}_{aa}} \!\!\!\!{(\left\|{x}_i-{x}_j\right\|^2-(d_{i j}^{\star 2}+\mu_{i j}))^2}+\|\Lambda(\boldsymbol{x}-(\boldsymbol{x}^{\star}+\boldsymbol{\epsilon}))\|^2
\end{align}

On this basis, we can establish that when the measurement errors are constrained, the local NE  $\hat{\boldsymbol{x}}$ is also the unique global NE  $\boldsymbol{x}^{\Diamond}$, which corresponds to the unique solution $\boldsymbol{x}^{\ddagger}$ of noisy SNL problem \eqref{noisySNL}. Moreover, the position errors between  NE $\boldsymbol{x}^{\Diamond}$ and actual positions $\boldsymbol{x}^{\star}$ can be quantified, as detailed in the following theorem, whose
 proof is given in Appendix \ref{a4n}.
% From the entire network’s perspective, it is desirable to obtain the optimum location selection Ψ
% It is desirable to obtain the optimum location selection $\boldsymbol{x}^{\star}$ from a global perspective. 
%In view of this, 
 \begin{mythm}\label{t11}
{
%  Consider fixed $\theta_{il}\in [0,\pi]$ for $l\in\mathcal{N}_a$, $\boldsymbol{\mu}$, $\boldsymbol{e}$, and $\boldsymbol{\epsilon}$.  
  Under Assumption 1, there exists a small positive $\delta$ and positive constants $a$,  $b$ and $c$ depending on $\delta$ such that  if 
 $\boldsymbol{\mu}$, $\boldsymbol{e}$, and $\boldsymbol{\epsilon}$ satisfy $a\|\boldsymbol{\mu}\|^2+b\|\boldsymbol{e}\|^2+c\|\boldsymbol{\epsilon}\|^2< \delta$, 
  %there exist    positive constants $a$,  $b$ and $c$ such that for any
% suitably small positive $\delta$, if 
 %Consider fixed 
 %$\theta_{il}\in [0,\pi]$ for $l\in\mathcal{N}_a$, 
 %$\boldsymbol{\mu}$, $\boldsymbol{e}$, and $\boldsymbol{\epsilon}$.
% There exist    positive constants $a$,  $b$, $c$ and $d$ 
 %such that for any
%suitably small positive $\delta$, if  $a\|\boldsymbol{\mu}_1\|^2+b\|\boldsymbol{\epsilon}\|^2+c\|{\boldsymbol{e}}_2\|^2+< \delta$,
%$\delta_1$, $\delta_2$, $\delta_3$ 
%and  positive constants $a$,  $b$ and $c$ depending on $\delta$
%and  positive constants $c$,  $d$ and $e$ depending on $\delta_1$, $\delta_2$ and $\delta_3$
%such that for any fixed $\boldsymbol{\mu}_1$, $\boldsymbol{\mu}_2$, and $\boldsymbol{\epsilon}$ satisfying 
% $a\|\boldsymbol{\mu}_1\|^2+b\|\boldsymbol{\mu}_2\|^2+c\|\boldsymbol{\epsilon}\|^2< \delta$,
%There exists a suitably small positive $\delta_1$, $\delta_2$, $\delta_3$ and  given fixed $\theta_{il}\in [0,\pi]$ for $l\in\mathcal{N}_a$
%an associated positive constant $c$ and $d$
%such that if the measurement errors in the squares of the distances obey $\|\boldsymbol{\mu}_1\|\leq \delta_1$, $\|\boldsymbol{\mu}_2\|\leq \delta_2$ and $\|\boldsymbol{\epsilon}\|\leq \delta_3$, 
then 
\begin{enumerate}[i)]
\item  there exists a unique NE $\boldsymbol{x}^{\Diamond}$, which is equal to   $\boldsymbol{x}^{\ddagger}$;
%in the noisy SNL problem \eqref{noisySNL};

\item  the NE $\boldsymbol{x}^{\Diamond}$ satisfies $\|\boldsymbol{x}^{\Diamond}-\boldsymbol{x}^{\star}\|^2< \delta$.
%${\Phi}(\boldsymbol{x}^{\Diamond})-{\Phi}(\boldsymbol{x}^{\star})< \delta$.
%$\|\boldsymbol{x}^{\Diamond}-\boldsymbol{x}^{\star}\|\rightarrow 0$ as $\|\boldsymbol{\mu}_1\|\rightarrow 0$, $\|\boldsymbol{\mu}_2\|\rightarrow 0$ and $\|\boldsymbol{\epsilon}\|\rightarrow 0$.
\end{enumerate}
}
%$\|\boldsymbol{x}^{\Diamond}-\boldsymbol{x}^{\star}\|\leq c\|\boldsymbol{\mu}\|+d\|\boldsymbol{\epsilon}\| $
\end{mythm}
 %\cite{calafiore2010distributed}
%\noindent\textbf{Proof sketch}
 
 %\subsection{Design scheme}
% 
%  Both $J_{i}$ and $	H$ 
% %	
% %we employ  the Least-Squares goodness of fit criterion  to measure the difference between the computed squared distance  and  the real squared distance,  
%

Theorem \ref{t11} shows that when the measurement noises are not large, the NE is unique and returns sensor position estimates that are close to the actual positions. Also, the errors between the actual positions and NE go to zero asymptotically as the noise perturbations in the actual positions and anchor positions go to zero. This establishes that a network can be approximately localized when the inter-node distance measurements and anchor node positions are contaminated with sufficiently small errors.

Regarding the error bound $\delta$, it is important to limit its magnitude. We propose the following algorithm to determine the value of $\delta$. Note that 
%$\Phi_1$, $\Phi_2$, $\delta_1$, $\delta_2$ and 
$R$ denotes a large enough constant introduced in the proof of Theorem \ref{t11} in Appendix \ref{a4n}. 
%The mini

%Based on Theorem \ref{t11}, we propose the following algorithm  to determine the value of $\delta$. 
%we put together all the procedures in Algorithm 1 to determine the value of $\delta$.

\begin{algorithm}[h]
	\renewcommand{\thealgorithm}{1}
	%\SetAlgoRefName{} % no count number
	\caption{}
	\label{a1}
	% \vspace{0.1cm}
%	\begin{algorithmic}
	\textbf{Initialization}:
 error bounds $\delta_1$, $\delta_2$, constants $a$, $b$, $c$ and $R$, 
 
 \quad\;  uncertainty vectors $\boldsymbol{\mu}$, $\boldsymbol{e}$, and $\boldsymbol{\epsilon}$ %with $a\|\boldsymbol{\mu}\|^2+b\|\boldsymbol{e}\|^2+$
 
% \quad\;
% %\quad\,\quad\, \;\,\;\;
% $c\|\boldsymbol{\epsilon}\|^2\leq \delta_1$ 
%and $a\|\boldsymbol{\mu}\|^2+b\|\boldsymbol{e}\|^2+c\|\boldsymbol{\epsilon}\|^2\leq \delta_2$
	% \vspace{-0.2cm}
 
%	aggregator $\sigma(\bm x)$, payoff functions $J_i(x_i,\sigma(\bm x))$, coupling constraint $\Omega_{\bm \alpha}$, all uncertain feasibility $\bm{\mathrm {A}}$, and a well-developed vGNE-seeking solver.
	% \vspace{-0.6cm}
	
	%\textbf{1)  }
 
%\quad\textbf{for} $k=1,2,\dots,M$, \textbf{do}

\textbf{1) }\;\textit{Input}:\quad $\delta_1$, $a$, $b$, and $c$, $\boldsymbol{\mu}$, $\boldsymbol{e}$ and $\boldsymbol{\epsilon}$ with $a\|\boldsymbol{\mu}\|^2+b\|\boldsymbol{e}\|^2+$

% \quad\,\quad\,\quad\, \;\,\;$\boldsymbol{\mu}$, $\boldsymbol{e}$ and $\boldsymbol{\epsilon}$ with $a\|\boldsymbol{\mu}\|^2+b\|\boldsymbol{e}\|^2+$

\quad\,\quad\,\quad\, \;\,\; $c\|\boldsymbol{\epsilon}\|^2\leq \delta_1$
%constants $a$, $b$, $c$ 
%uncertain $\alpha_i[k]\in\mathrm{A}_i,~i\in\mathcal I$,

\quad\, \textit{Solve}:\quad
%problem \eqref{p11}
%\quad\,\quad\,\quad\, \;\,\; 
\textbf{while} $\nabla^2 {\Phi}(\boldsymbol{x}, \boldsymbol{\mu}, \boldsymbol{e}, \boldsymbol{\epsilon})\preceq 0 $ \textbf{do}
%for  $a\|\boldsymbol{\mu}\|^2+b\|\boldsymbol{e}\|^2+$

% \quad\,\quad\,\quad\, \;\,\;\quad\, $c\|\boldsymbol{\epsilon}\|^2\leq \delta_1$} \textbf{do}

\quad\,\quad\,\quad\, \;\,\;\quad\, reset $ \delta_1$,  $\boldsymbol{\mu}$, $\boldsymbol{e}$ and $\boldsymbol{\epsilon}$ such that $a\|\boldsymbol{\mu}\|^2+b\|\boldsymbol{e}\|^2+$

\quad\,\quad\,\quad\, \;\,\;\quad\, $c\|\boldsymbol{\epsilon}\|^2< \delta_1$

\quad\,\quad\,\quad\, \;\,\;  \textbf{end while} 

%\quad\,\quad\,\quad\, \;\,\;\quad\, {choose $\delta_2 $ such that $\Phi_2= \frac{1}{2}\Phi_1$}

\quad\,\quad\,\quad\, \;\,\;  solve  \eqref{p11} and obtain $\Phi_1$ and $\delta_1$
% $
% %\boldsymbol{x}^{\prime}=
% {\operatorname{min}}_{{\|\boldsymbol{x} -\boldsymbol{x}^{\star}\|^2\geq \delta_1}}
% %\in \mathcal{B}^c} 
% {\Phi}_1(\boldsymbol{x} )$
%\quad\quad\quad\quad\;\,\quad seeking solver,

\quad\, \textit{Output}: $\frac{1}{2}\Phi_1$, $\delta_1$
%vGNE $\bm x^*_{\bm \alpha[k]}$ of game $\mathscr G_{\bm \alpha[k]} (\boldsymbol{x}^{\prime})$.

%\quad\textbf{end for}

	%\textbf{2) Inverse learning}
 
\textbf{2) }\;\textit{Input}:\quad $\delta_1$, $\delta_2$, $a$, $b$, $c$, $R$, $\frac{1}{2}\Phi_1$

%\quad\,\quad\,\quad\, \;\,\; measurement errors $\boldsymbol{\mu}$, $\boldsymbol{e}$ and  $\boldsymbol{\epsilon}$, $\frac{1}{2}\Phi_1$
%with $a\|\boldsymbol{\mu}\|^2+$

%\quad\,\quad\,\quad\, \;\,\;\;$b\|\boldsymbol{e}\|^2+c\|\boldsymbol{\epsilon}\|^2\leq \delta_2$, $\frac{1}{2}\Phi_1$
%all data $(\bm \alpha[k],\bm x^*_{\bm \alpha[k]}),~k=1,2,\dots,M$,

\quad\, \textit{Solve}:\quad 
\textbf{while} $\delta_2>0$ \textbf{do}

% $
% (\boldsymbol{x}^{\prime \prime }, \boldsymbol{\mu}^{\prime}, \boldsymbol{\epsilon}^{\prime}, \boldsymbol{e}^{\prime})=$
%\quad\,\quad\,\quad\,
\quad\,\quad\,\quad\, \;\,\;\quad\,  solve  \eqref{phi2} and obtain ${\Phi}_2$
% $
% {\operatorname{min}}_{{\delta_1 \leq\|\boldsymbol{x} -\boldsymbol{x}^{\star}\|^2\leq R, a\|\boldsymbol{\mu}\|^2+b\|\boldsymbol{e}\|^2+c\|\boldsymbol{\epsilon}\|^2\leq \delta_2
% %\boldsymbol{x} \in \mathcal{B}^c, 
% }} 
% %{\Phi}_2(\boldsymbol{x})$
% %\nolimits_{} 
% \;{\Phi}_2(\boldsymbol{x},\boldsymbol{\mu},\boldsymbol{\epsilon}, \boldsymbol{e})$
%inverse VI-based  learning approach in \eqref{inverse_data},

\quad\,\quad\,\quad\, \;\,\;\quad\, \textbf{if} {$\Phi_2\geq \frac{1}{2}\Phi_1$ for all $\delta_2 \leq \delta_1$} 

\quad\,\quad\,\quad\, \;\,\;\quad\,\quad  $\delta_2 := \delta_1$

\quad\,\quad\,\quad\, \;\,\;\quad\,\quad \textbf{break} 

\quad\,\quad\,\quad\, \;\,\;\quad\,   \textbf{else} 

\quad\,\quad\,\quad\, \;\,\;\quad\,\quad  {choose $\delta_2 $ such that $\Phi_2= \frac{1}{2}\Phi_1$}

\quad\,\quad\,\quad\, \;\,\;\quad\,\quad \textbf{break}

\quad\,\quad\,\quad\, \;\,\;\quad\, \textbf{end if} 

\quad\,\quad\,\quad\, \;\,\;  \textbf{end while} 

    \quad \;\,\textit{Output}: \,$\delta_2$
    %estimated weight $\hat{\bm \beta}$ in black-box aggregator

	%\textbf{3) Robust Counterpart}
 
\textbf{3) }\,\,\textit{Input}:\; \;\,error bound $\delta_2$, $\frac{1}{2}\Phi_1$
%\quad all uncertain flexibility $\bm{\mathcal U}$ and estimated weight $\hat{\bm \beta}$ 

% \quad\quad\quad\;\,\quad within the aggregator,

  \quad \;\,\textit{Solve}:\; \;\,\,$\delta:= \min\{\delta_2,\frac{1}{2}\Phi_1\}$
%the deterministic game $\hat{\mathscr  G}$ in \eqref{robust_counter} with auxiliary 

%\quad\quad\quad\;\,\quad variables $y_i,i\in\mathcal  I$ and estimated weight $\hat{\bm \beta}$,

% \quad\quad\quad\;\,\quad $y_i,i\in\mathcal  I$,

  \quad \;\,\textit{Output}:  \;$\delta$
%r-GNE $\bm x^*$ of game $\mathscr G$.

%\end{algorithmic}
\end{algorithm}

 There are three reasons for limiting the value of $\delta$.
 %in Algorithm 1:
%As for the error bound $\delta$, we need to restrain its value for three reasons.

%As for the value of $\delta$ there are three reasons.
%for limiting this quantity. 

1) To ensure the positive definiteness of $\nabla^2 \Phi$.
%so that the implicit function theorem is valid.
%which means that $\nabla \Phi$ should be positive definite.
%This property is satisfied at the point 
Note that $\nabla^2 \Phi$ is positive definite at $(\boldsymbol{\mu}_0,\boldsymbol{e}_0,\boldsymbol{\epsilon}_0)$, but may cease to hold well away from this point. Thus, we need to limit the size of $\delta_1$ in Step 1 of Algorithm 1. 

2) To avoid getting stuck into local NE in noisy sensor localization.
%under noisy conditions.
In the noiseless case, there may exist a   local NE (denoted as $\hat{\boldsymbol{x}}$), whose value  $\Phi(\hat{\boldsymbol{x}})$ differs from zero by a small amount. 
The coordinate values corresponding to local NE $\hat{\boldsymbol{x}}$ are different from those corresponding to the global NE ${\boldsymbol{x}}^{\Diamond}={\boldsymbol{x}}^{\star}$. Then, when the noise levels are slowly increased from zero, the coordinate values corresponding to the global NE ${\boldsymbol{x}}^{\Diamond}$  may jump at some noise level to the coordinate values corresponding to local NE $\hat{\boldsymbol{x}}$. To avoid this situation, we should take $a\|\boldsymbol{\mu}\|^2+b\|\boldsymbol{e}\|^2+c\|\boldsymbol{\epsilon}\|^2\leq \delta_2< \delta_1$  to ensure $\Phi(\hat{\boldsymbol{x}})$  is not smaller than $\frac{1}{2}\Phi_1$. Therefore,  $\Phi(\hat{\boldsymbol{x}})$ should be larger enough than  $\Phi({\boldsymbol{x}}^{\Diamond})$. 

3)  To ensure that the global NE $\boldsymbol{x}^{\Diamond}$  in  noisy sensor localization is close to the real localization ${\boldsymbol{x}}^{\star}$. To achieve this, we should take $a\|\boldsymbol{\mu}\|^2+b\|\boldsymbol{e}\|^2+c\|\boldsymbol{\epsilon}\|^2< \frac{1}{2}\Phi_1$,   so as to let  $\Phi({\boldsymbol{x}}^{\Diamond})$  not greater than $\frac{1}{2}\Phi_1$. Hence, $\boldsymbol{x}^{\Diamond}$ is not far from the actual positions ${\boldsymbol{x}}^{\star}$.

%compute an r-GNE within a black-box aggregative game problem (2)

\vspace{-0.2cm}
\begin{remark}
%As for the value of $a$, $b$ and $c$ 
To determine the values of $a$, $b$ and $c$, we can follow two steps. Firstly, we can employ the implicit function theorem or Lipschitz continuity to obtain $a$, $b$, and $c$ based on Lemma \ref{ee1}. Secondly, we compare them respectively with $a^{\prime}=1$, $b^{\prime}=2$ and $c^{\prime}=32d_{max}^{\star2}$ (see the last line of (\ref{pp3}) in Appendix \ref{a4n}). If $a<1$ or $b<2$ or $c<32d_{max}^{\star2}$, then we set $a=1$, $b=2$ and $c=32d_{max}^{\star2}$; otherwise, we keep them unchanged. The reason for limiting  $a$, $b$, and $c$ is not only to satisfy $\|\boldsymbol{x}^{\Diamond}-\boldsymbol{x}^{\star}\|<\delta$ but also to ensure  $\Phi({\boldsymbol{x}}^{\Diamond})<\frac{1}{2}\Phi_1$.
\end{remark}

Additionally, in the case that the anchor's positions are perfectly known, i.e., $\boldsymbol{\epsilon}=\boldsymbol{\epsilon}_0=\operatorname{col}\{0\}_{l=1}^{|\mathcal{N}_a|}$, the results in Theorem \ref{t11} can be further simplified.
\begin{corollary}\label{c1}
 Under Assumption 1, there exists a small positive $\delta$ and positive constants $\bar{a}$ depending on $\delta$ such that  if $\boldsymbol{\epsilon}=\boldsymbol{\epsilon}_0$ and 
 $\boldsymbol{\mu}$ and $\boldsymbol{e}$ satisfy $\|\boldsymbol{\mu}\|+\|\boldsymbol{e}\|< \delta$, 
then the NE of SNL problem \eqref{f1} is
unique and $\|\boldsymbol{x}^{\Diamond}-\boldsymbol{x}^{\star}\|\leq \bar{a}(\|\boldsymbol{\mu}\|+\|\boldsymbol{e}\|)$.
\end{corollary}	
% \textbf{Proof}
Also, when the inter-node distance measurement errors are zero, i.e., $\boldsymbol{\mu}=\boldsymbol{\mu}_0=\operatorname{col}\{\mu_{ij}=0\}_{(i,j)\in \mathcal{E}_{ss}}$ and $\boldsymbol{e}=\boldsymbol{e}_0=\operatorname{col}\{e_{il}=0\}_{(i,l)\in \mathcal{E}_{as}}$, we have the following results.
\begin{corollary}\label{c2}
 Under Assumption 1, there exists a small positive $\delta$ and positive constants $\bar{b}$ depending on $\delta$ such that  if $\boldsymbol{\mu}=\boldsymbol{\mu}_0$, $\boldsymbol{e}=\boldsymbol{e}_0$ 
 and 
 $\boldsymbol{\epsilon}$ satisfies $\|\boldsymbol{\epsilon}\|< \delta$, 
then the NE of SNL problem \eqref{f1} is
unique and  $\|\boldsymbol{x}^{\Diamond}-\boldsymbol{x}^{\star}\|\leq \bar{b}\|\boldsymbol{\epsilon}\|$.
\end{corollary}

\section{Numerical Experiments}
	\begin{figure}[ht]
			\hspace{-0.9cm}
			\centering	
			\subfigure[True localization]{
				\begin{minipage}[t]{0.45\linewidth}
					\centering
					\includegraphics[width=4.6cm]{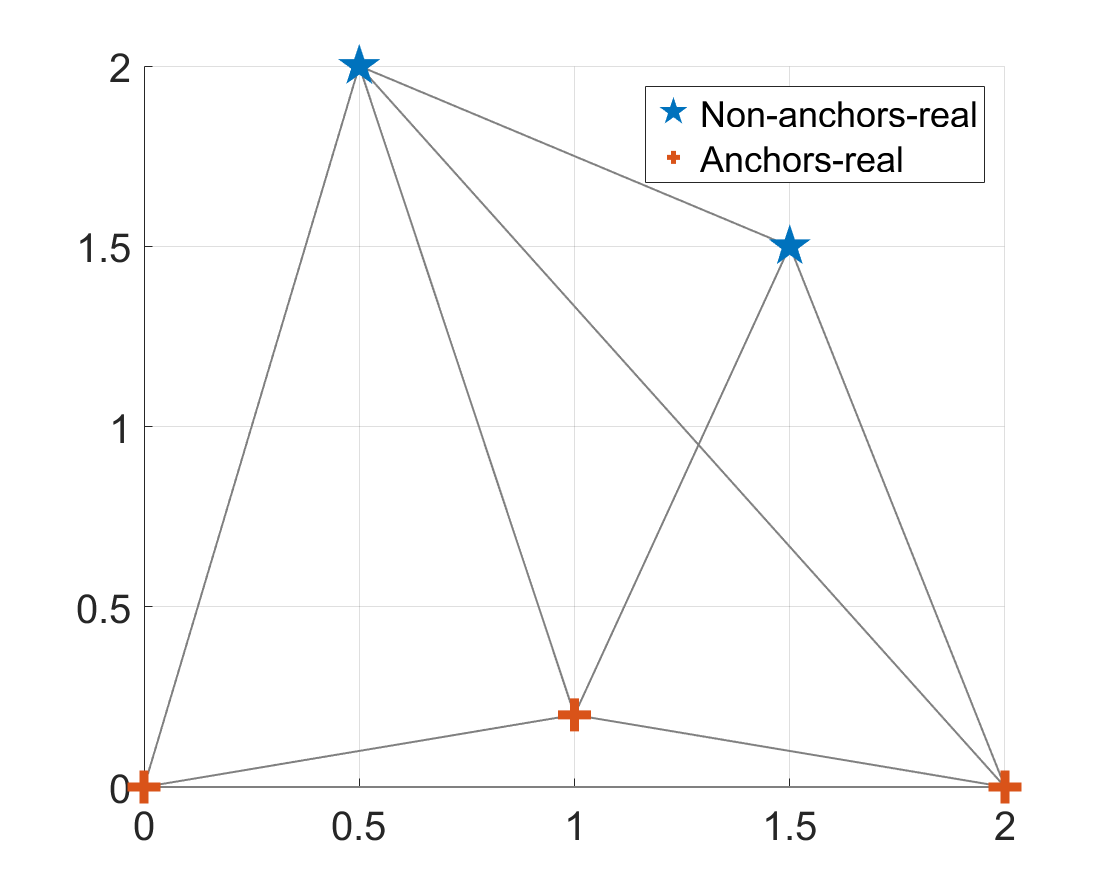}
					%\caption{fig1}
				\end{minipage}%
			}%
			\hspace{0.1cm}
			\subfigure[Noisy localization]{
				\begin{minipage}[t]{0.45\linewidth}
					\centering
					\includegraphics[width=4.6cm]{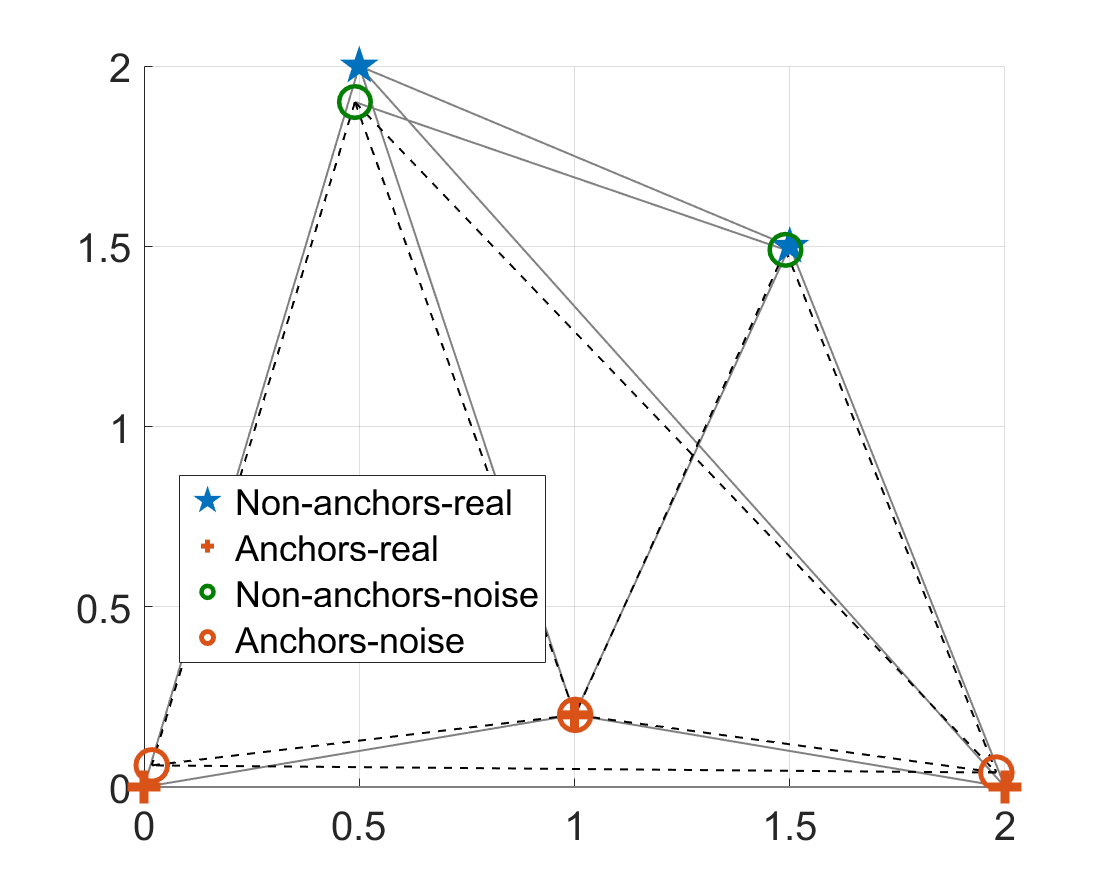}
					%\caption{fig2}
				\end{minipage}%
			}%
			\centering
   \vspace{-0.1cm}
			\caption{The true localization and noisy localization results.} 
   \vspace{-0.3cm}
			\label{fig12345s}
		\end{figure}
    \begin{figure*}[hb]
			%\hspace{-0.7cm}
			\centering	
			\subfigure[$\delta=0.01$]{
				\begin{minipage}[t]{0.248\linewidth}
					\centering
					\includegraphics[width=4.7cm]{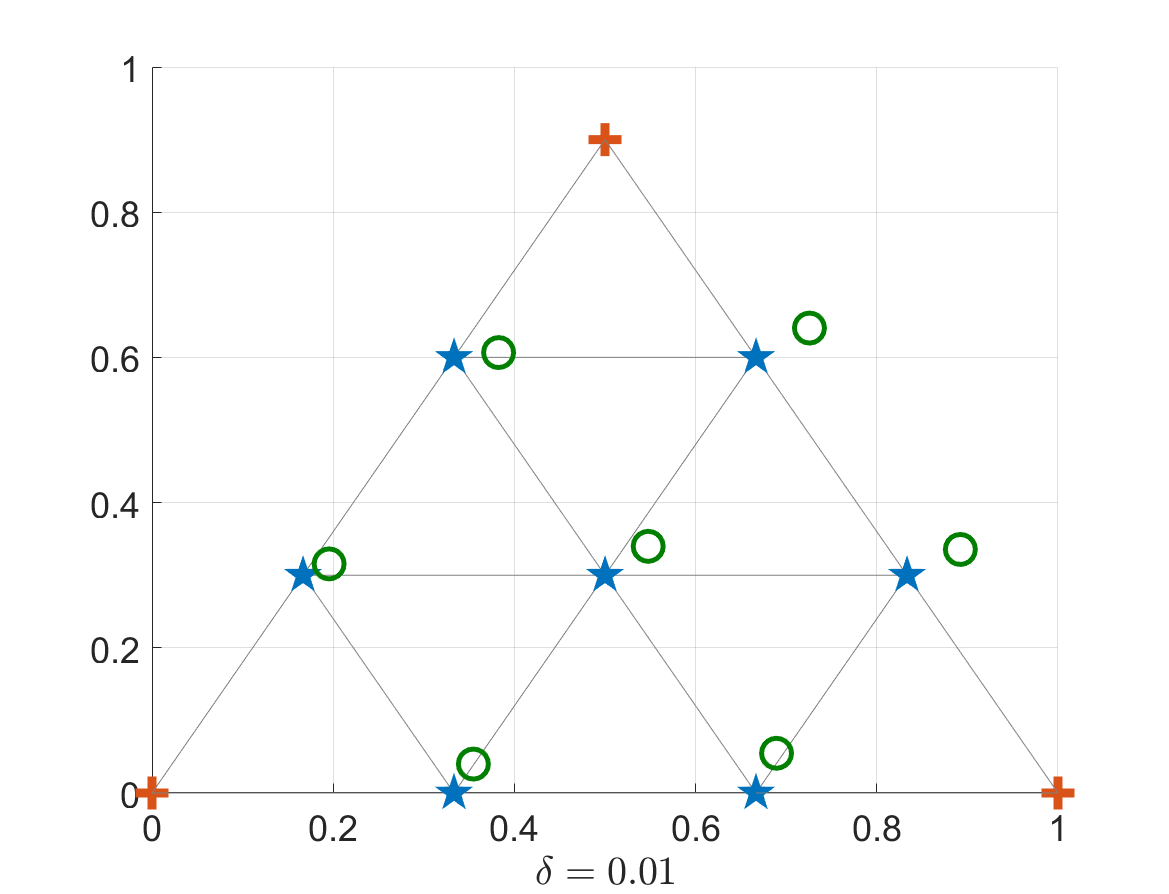}
    % \vspace{0.2cm}
					%\caption{fig1}
				\end{minipage}%
			}%
			%\hspace{0.2cm}
			\subfigure[$\delta=0.05$]{
				\begin{minipage}[t]{0.248\linewidth}
					\centering
					\includegraphics[width=4.7cm]{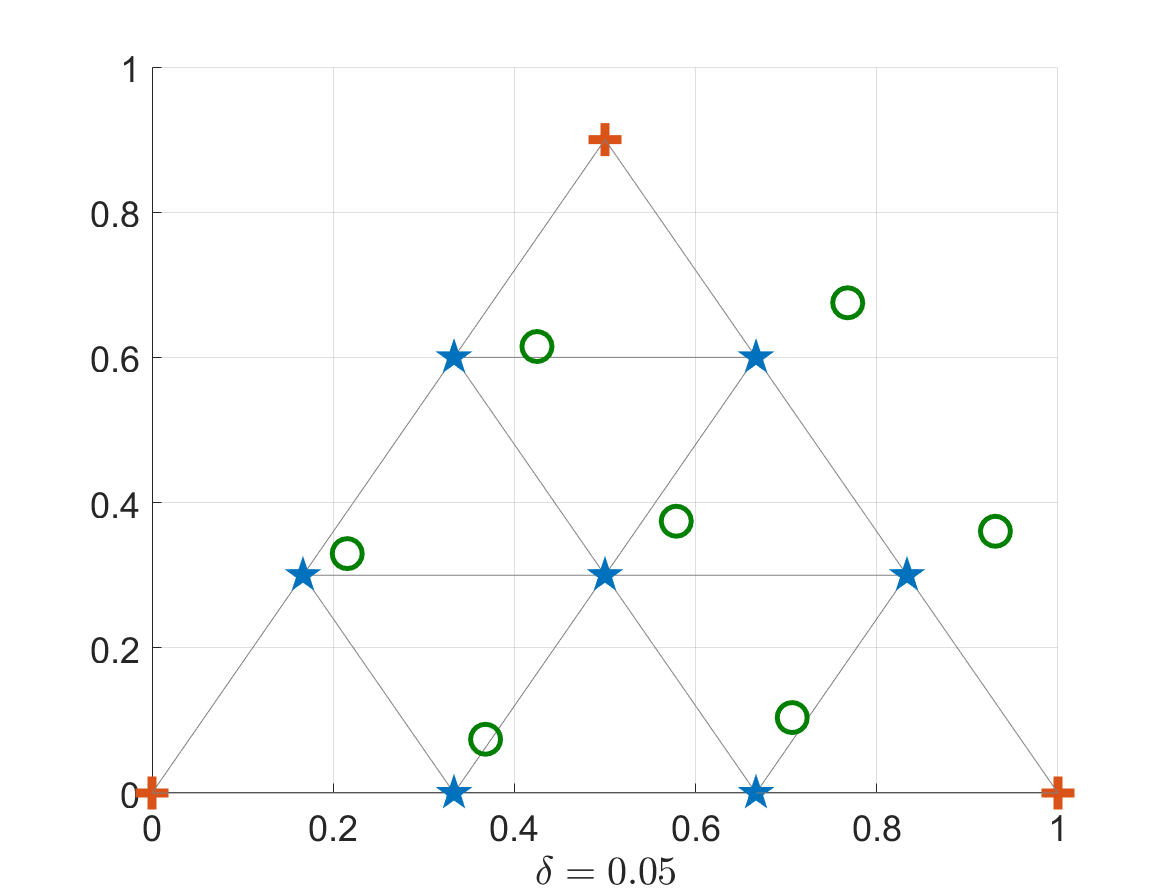}
   %  \vspace{-0.3cm}
					%\caption{fig2}
				\end{minipage}%
			}%
%   \\
   \subfigure[$\delta=0.1$]{
				\begin{minipage}[t]{0.248\linewidth}
					\centering
					\includegraphics[width=4.7cm]{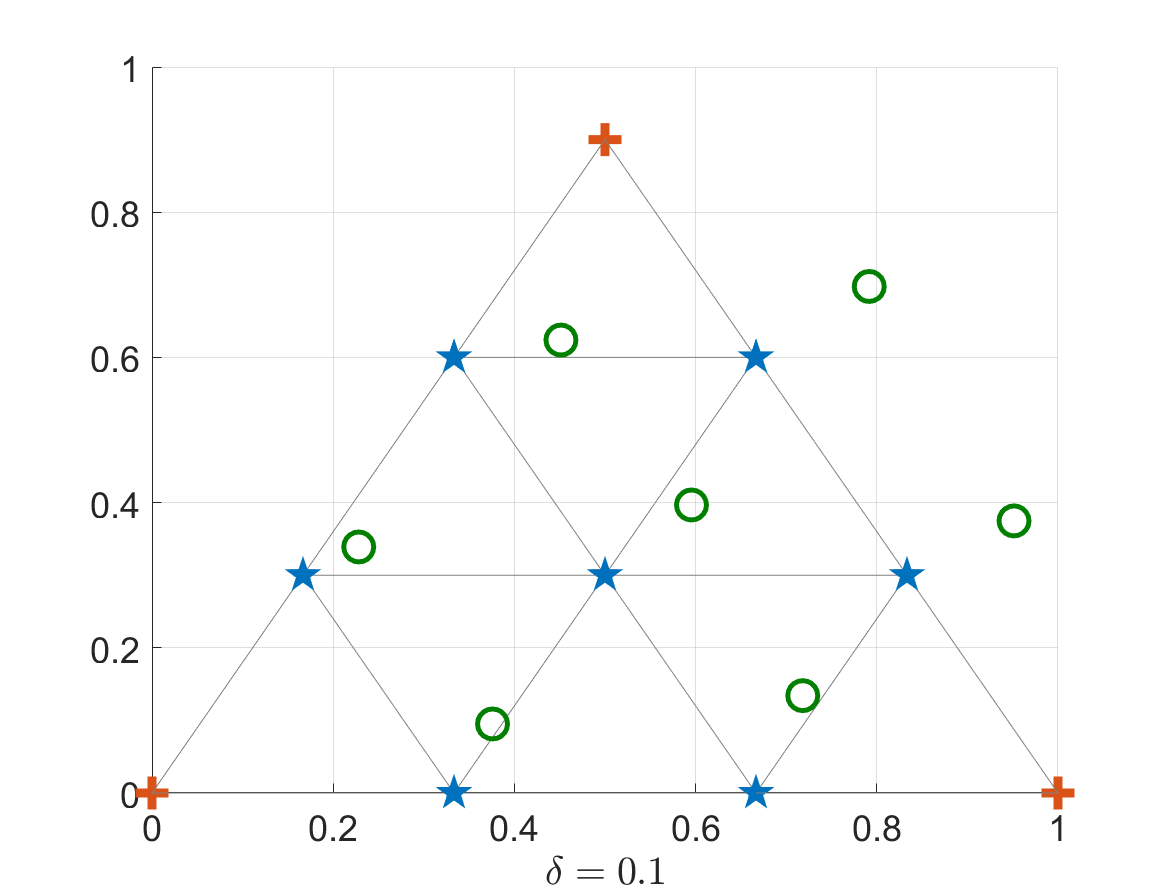}
     \vspace{-0.2cm}
					%\caption{fig1}
				\end{minipage}%
			}%
	%		\hspace{0.2cm}
			\subfigure[$\delta=0.2$]{
				\begin{minipage}[t]{0.248\linewidth}
					\centering
					\includegraphics[width=4.7cm]{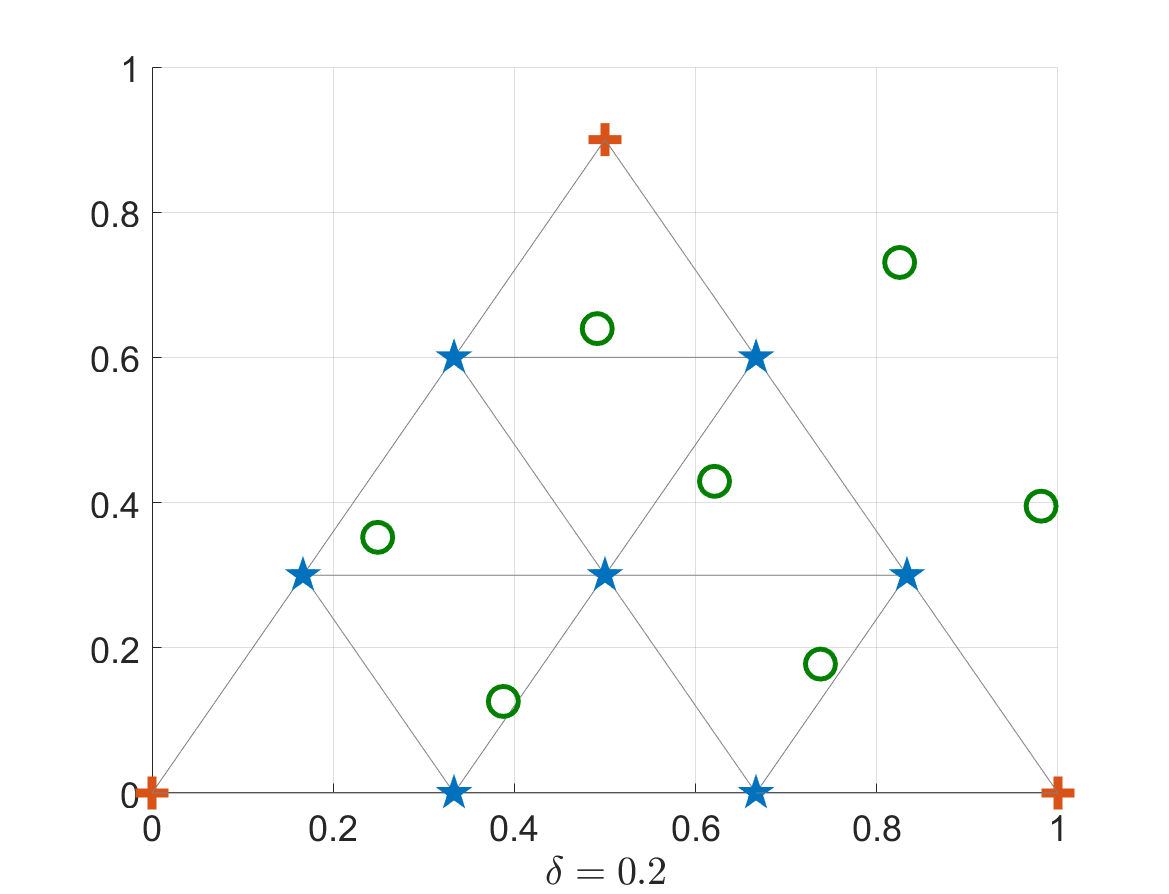}
					%\caption{fig2}
				\end{minipage}%
			}\\
   	\subfigure[$\delta=0.5$]{
				\begin{minipage}[t]{0.248\linewidth}
					\centering
					\includegraphics[width=4.7cm]{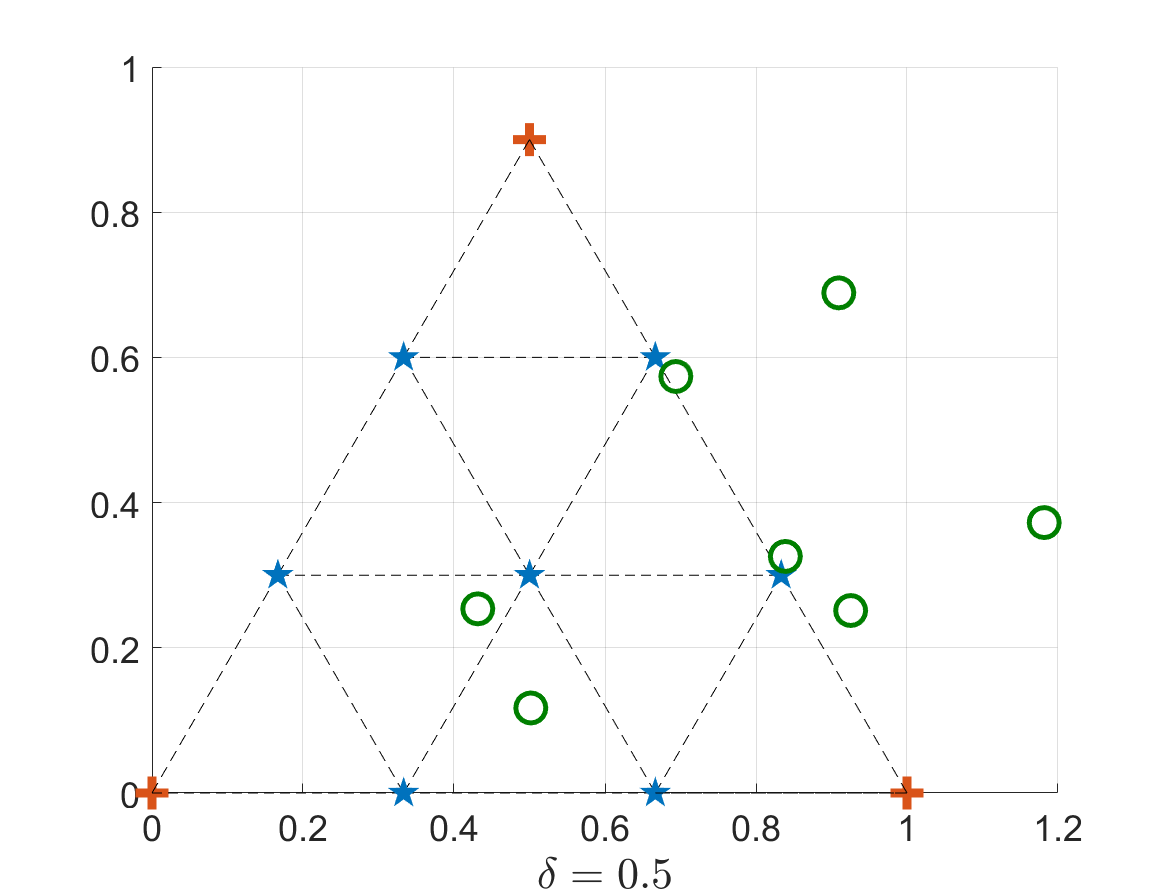}
    % \vspace{0.2cm}
					%\caption{fig1}
				\end{minipage}%
			}%
			%\hspace{0.2cm}
			\subfigure[$\delta=1$]{
				\begin{minipage}[t]{0.248\linewidth}
					\centering
					\includegraphics[width=4.7cm]{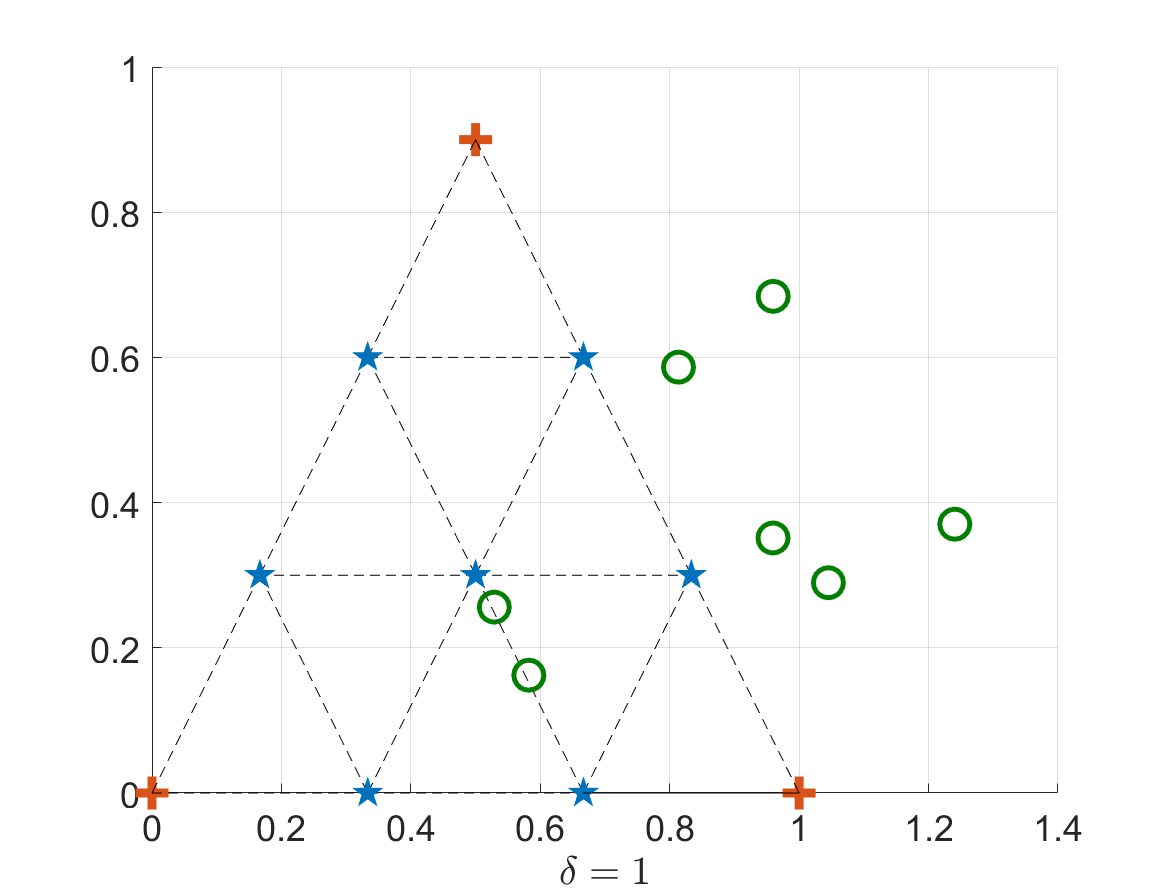}
   %  \vspace{-0.3cm}
					%\caption{fig2}
				\end{minipage}%
			}%
%   \\
   \subfigure[$\delta=1.5$]{
				\begin{minipage}[t]{0.248\linewidth}
					\centering
					\includegraphics[width=4.7cm]{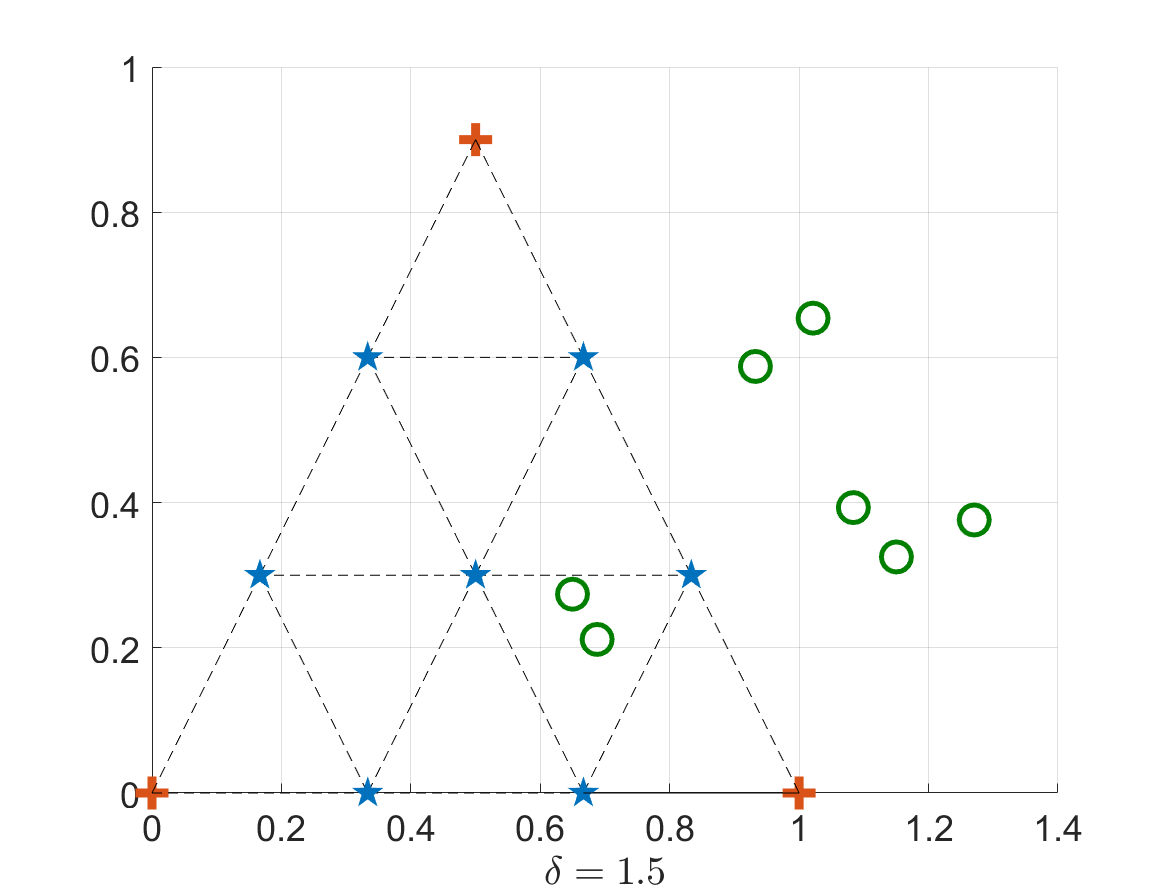}
     \vspace{-0.2cm}
					%\caption{fig1}
				\end{minipage}%
			}%
	%		\hspace{0.2cm}
			\subfigure[$\delta=2$]{
				\begin{minipage}[t]{0.248\linewidth}
					\centering
					\includegraphics[width=4.7cm]{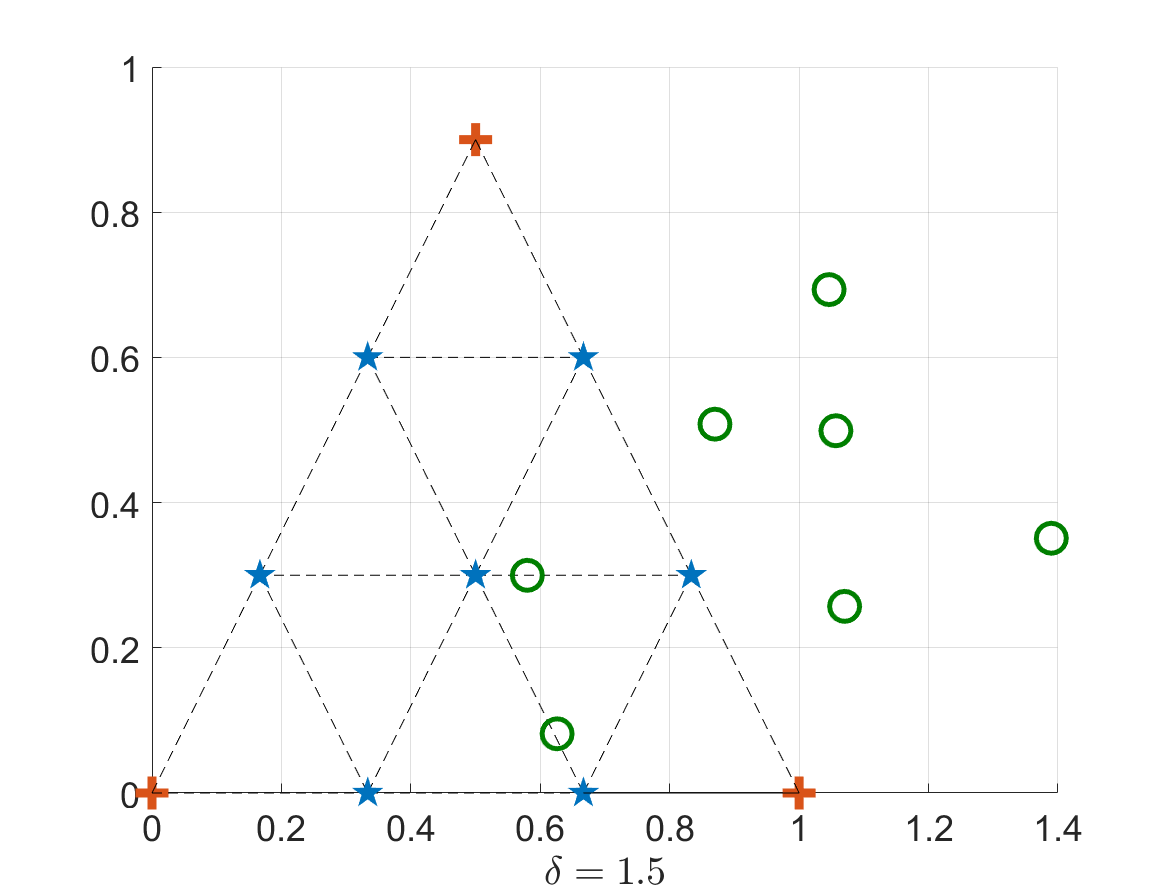}
					%\caption{fig2}
				\end{minipage}%
			}
			\centering
   \vspace{-0.2cm}
			\caption{Computed sensor location results with different error levels.} 
			\label{figgg4}
		\end{figure*}

\begin{table*}[t]
	\centering
	%	\scriptsize
	%	\small
	\footnotesize
	\caption{MLE of different measurement noise bounds.}
 \vspace{-0.1cm}
		\setlength\tabcolsep{20pt}
		\renewcommand\arraystretch{1.5}
		%	 \resizebox{\textwidth}{7mm}
		\begin{tabular}{c|c|c|c|c|c}
			\hline
			\hline
			\specialrule{0em}{0.2pt}{0.5pt}
	% 		\multirow{2}{1.5cm}{Numbers} & 	\multirow{1}{1.8cm}{$ N=10$}&  \multirow{1}{1.8cm}{$ N=50$}&  \multirow{1}{1.8cm}{$ N=80$} & \multirow{1}{1.8cm}{$ N=100$}\\ 
 % {} & 	\multirow{1}{1.8cm}{$ M=40$}&  \multirow{1}{1.8cm}{$  M=50$}& \multirow{1}{1.8cm}{$ M=100$} &\multirow{1}{1.8cm}{$ M=100$}\\
  &$ \delta=0.05 $ &$ \delta=0.1 $  &$  \delta=0.2 $& $  \delta=0.5 $ & $ \delta=1$ \\ 
   \hline	
			%		Algorithm 1 &0.05 &0.11 &1.04 &2.45  &5.14 \\ \hline
			%		Algorithm-SQP&0.17&0.36 &3.09 & 8.31& 17.78 \\  \hline
			%		Algorithm-IPM&0.36 &0.64 &5.94 &10.51& 24.83 \\\hline
			MLE&  $8.4433*10^{-4}$ &$4.0300*10^{-4}$  & 0.0262& 0.0561 &0.1468  \\ \hline
 % &     &   & &  \\ \hline
% ${\Phi}(\boldsymbol{x}^{\Diamond})-{\Phi}(\boldsymbol{x}^{\star})$& 0.0021 & 0.004  & & &  \\ \hline
			% PGD&    2.2700&2.2052  & 0.1891& 0.2419\\  \hline
			% penalty&2.2696 &2.7456  &0.0817 &0.3294 \\\hline
   % proximal&3.2578 &2.3697 &0.6794 &0.4881 \\\hline
			\hline
			%		Running time (sec)& 1.7747& 1.7358& 1.6012
			%		& 2.1808& 1.6239
			%		&  2.7643\\	\hline	
	\end{tabular}
	\label{tab6435}
 \vspace{-0.3cm}
\end{table*}
Consider an SNL problem with noisy measurements. Take $N=2$ non-anchor nodes and $M=3$ anchor nodes. The original configuration is shown in Fig. \ref{fig12345s}(a). The true positions of anchor and non-anchor nodes are represented by blue stars and red pluses respectively, consistent with Fig. \ref{fig123s}. The gray lines indicate connections between nodes. 
Following the procedures in Alg.1, we set $\delta_1=0.1$, $a=1$, $b=32$, $c=2$
and get $\delta=0.1$. The corresponding noisy localization result is shown in Fig. \ref{fig12345s}(b). 
%We can see that  when the measurement errors satisfy the relation in Theorem \ref{t11},
We can see from Fig. \ref{fig12345s}(b) that when the measurement errors satisfy the condition in Theorem \ref{t11}, the position errors between computed results and actual ones are also small. 
Moreover, Table \ref{tab6435} lists the values of \textit{mean localization error} (MLE) when measurement errors are constrained by different bounds, where  
\begin{equation*}
  \operatorname{MLE}=\frac{1}{N} \sqrt{\sum\nolimits_{i=1}^{N}\|x_{i}^{\Diamond}-x_{i}^{\star}\|^2}.
\end{equation*}
We can observe from Table \ref{tab6435} that with the increasing bound $\delta$, the MLE value also increases. Additionally, the MLE is consistently smaller than $\delta$.

Next, we consider another configuration with $N = 7$ non-anchor nodes and $M=3$ anchor nodes \cite{xu2024global}. Fig. \ref{figgg4} shows the noisy localization results under different error bounds. As $\delta$ ranges from 0.01 to 0.2, the noisy localization result (shown as green circles) gradually deviates from the true values (shown as blue stars), while its graph structure remains globally rigid, similar to the true one. However, when $\delta$ exceeds 0.5, the graph structure gradually becomes distorted. Particularly, at $\delta=2$, such measurement errors completely disrupt the final computation results.

%We use the same symbols in  Fig. \ref{fig12345s}.
%Moreover, 
%
% comparing
% two indexes, root mean square error
% The
% performance of this game setting can be evaluated by 
%the \textit{mean localization error}:

%, while the value of $\Delta$ goes rapidly (exponentially) to zero. 

%  \begin{figure}[ht]
% 	\centering	
% 	\includegraphics[scale=0.23]{n6.png}\\
%  \vspace{-0.2cm}
% 	\caption{Noisy SNL.}
% 	\label{fig12345s}
% 	\vspace{-0.3cm}
% \end{figure}

% Different network configurations
  
% The relation between different $\delta$ and $\|{\Phi}(\boldsymbol{x}^{\Diamond})-{\Phi}(\boldsymbol{x}^{\star})\|$

%xxxxxxxxxxxxxxxxxxxxxxxxxxxxxxxxx

\section{Concluding Remarks}
In this paper, we have studied the solution of the SNL problems with noisy measurements.  By formulating a non-convex SNL potential game, we have investigated the existence
and uniqueness of NE in both 
ideal settings with accurate anchor node location information and accurate inter-node distance measurements and practical settings
with anchor location inaccuracies and distance measurement noise.  In the noiseless case, we  have shown that the NE exists and
is unique, and corresponds to the precise network localization. In the case involving errors in anchor node position and inter-node distance measurements,  we have established that if these errors
are sufficiently small, the NE exists and is still unique, providing
an approximate solution to the SNL problem. Moreover, the position errors from  NE to the precise network localization can be quantified, providing that the measurement errors are constrained by a small bound.

% . Then based on the canonical duality theory, we have proposed  a conjugation-based algorithm  to compute the stationary point of a complementary dual problem, which  actually induces the global
% NE  if a duality relation can be checked. Finally, the computational efficiency
%  of our algorithm has been
% illustrated by several experiments. 

In follow-up works, we plan to extend our results to more complex scenarios including i) designing an algorithm with global convergence for  NE seeking, ii) generalizing the model to a distributed case, and iii) exploring milder graph conditions.
% ii) generalizing the model to cases
% with measurement noise, and 
\appendices
\vspace{-0.2cm}
\section{Graph Theory  }\label{lea1}
We first consider {an $n$-dimensional  representation of sensor network graph $\mathcal{G}$, which is a mapping of  $\mathcal G(\mathcal N, \mathcal E)$  to the point 
formations $\bar{\boldsymbol{x}}: \mathcal N \rightarrow \mathbb{ R}^{n}$, where %${\boldsymbol{x}}=\operatorname{col}\{x_1,\dots,x_{Q}\}$ and 
$\bar{\boldsymbol{x}}(i)=x_i^{\operatorname{T}}$ is  the row vector of the coordinates of
the $i$-th node in $\mathbb{ R}^{n}$ and $x_i\in \mathbb{ R}^{n}$. 
%In this paper, the  $x_i$ is the actual position of sensor node $i$.
%and  $ \mathcal{E}=\{(i,j)\in\mathcal{N}\times\mathcal{N}:\|x_{i}-x_{j}\|\leq R_s,i\neq j\} $ is the edge set between sensors, where $ R_s $ is a fixed range representing a sensor's capability of sensing range measurements from others.
% Also each sensor has the capability of sensing range measurements from other sensors within a fixed range $ R_s $. 
% This indicates that  
%Consider a map $\bar{p}: \mathcal N \rightarrow \mathbb{ R}^{n}$ with $n\in\{2,3\}$, which assigns a location in $\mathbb{ R}^{n}$ to each node in $\mathcal N$. 
Given the graph  $\mathcal G(\mathcal N, \mathcal E)$ and an $n$-dimensional representation $\bar{\boldsymbol{x}}$ of it, the pair $(\mathcal G, \bar{\boldsymbol{x}})$ is called a $n$-dimensional framework. A framework  $(\mathcal G, \bar{\boldsymbol{x}})$ is called generic\footnote{Some special configurations exist among the sensor positions, e.g., groups of sensors may be collinear. The reason for using the term generic is to highlight the need to exclude the problems arising from such configurations.} if the set containing the coordinates of all its points is algebraically independent over the
rationales \cite{anderson2010formal}.  {A framework $(\mathcal G, \bar{\boldsymbol{x}})$ is called  rigid if there exists a sufficiently small positive constant $\epsilon$ such that if every framework $(\mathcal G, \bar{\boldsymbol{y}})$
satisfies  $\|x_i-y_i \|\leq\epsilon$ for $i\in\mathcal N$ and 
$\| x_i-x_j\|=\|y_i-y_j\|$ for every  pair $i,j\in \mathcal N$ connected by an edge in $\mathcal E$, then 
$\| x_i-x_j\|=\|y_i-y_j\|$  for any node pair $i,j\in \mathcal N$  no matter there is an edge between them. Graph
$\mathcal G(\mathcal N, \mathcal E)$  is called generically   $n$-rigid or simply rigid (in $n$ dimensions) if any generic framework 
$(\mathcal G, \bar{\boldsymbol{x}})$ is rigid. A framework $(\mathcal G, \bar{\boldsymbol{x}})$  is globally rigid if every framework $(\mathcal G, \bar{\boldsymbol{y}})$ satisfying
$\| x_i-x_j\|=\|y_i-y_j\|$  for any node pair $i,j\in \mathcal N$  connected by an edge in $\mathcal E$ and  $\| x_i-x_j\|=\|y_i-y_j\|$  for any node pair $i,j\in \mathcal N$  that are not connected by a single edge. %between them.
Graph
$\mathcal G(\mathcal N, \mathcal E)$  is called generically globally rigid if 
%there is an associated framework
any generic framework 
$(\mathcal G, \bar{\boldsymbol{x}})$  is 
%generic and 
globally rigid \cite{anderson2010formal,fidan2010closing,tay1985generating}.}

%\vspace{-0.2cm}

\section{Proof of Theorem \ref{l11}}\label{rea2}
%\vspace{-0.2cm}
Recall the formulation of the non-convex function $\Phi$. 
 If all the non-anchor nodes are localized accurately such that $\|x_{i}^{\star}-x_{j}^{\star}\|^2-d_{i j}^2=0$ for any  $(i,j)\in \mathcal{E}$, then 
$\Phi(\boldsymbol{x}^{\star})$ tends to be zero. Thus  $\boldsymbol{x}^{\star}$ is the global minimum of $\Phi$. 
%The optimization problem to find the actual positions of non-anchor nodes can be regarded as finding the global minimum of $\Phi$. 
%Moreover, referring to \cite{eren2004rigidity, aspnes2006theory},  Assumption 1 guarantees that the global minimum of $\Phi$ is unique and corresponds to the actual position profile $\boldsymbol{x}^{\star}$. 
Moreover, referring to  \cite{eren2004rigidity, aspnes2006theory}, Assumption 1 guarantees that the global minimum of $\Phi$ is unique and corresponds to the actual position profile $\boldsymbol{x}^{\star}$. 
%This completes the proof
%On the other hand, 
Referring to Lemma 1, 
a global NE $\boldsymbol{x}^{\Diamond}$ is equal to a global minimum of the potential function $\Phi$. Since the global minimal of $\Phi$ corresponds to the uniquely actual localization $\boldsymbol{x}^{\star}$, we can obtain the conclusion. \hfill $\square$

%\vspace{-0.4cm}

\section{Proof of Lemma \ref{ee1}}\label{rea}

%\vspace{-0.2cm}

We need  additional notations for the proof of Lemma \!\ref{ee1}.

Define 
the rigidity matrix by $R$ with $|\mathcal{E}|$ rows and $2|\mathcal{N}|$ columns. Each edge gives rise to a row, and if the edge links vertices $j$ and $k$,
the nonzero entries of the row of the matrix are in columns $2j- 1$, $2j$, $2k-1$, and $2k$
and are, respectively, $x_{j1}-x_{k1}$, $x_{j2}-x_{k2}$, $x_{k1}-x_{j1}$, $x_{k2}-x_{j2}$.  Define the reduced
rigidity matrix $R_r$ to be the submatrix of $R$ containing those columns corresponding
to vertices $\mathcal{N}_s=1, 2, \dots, N$ and those edges joining vertex pairs of which at least one is in
the set ${1, 2, \dots, N}$. Suppose that the vertices are ordered so that the last numbered vertices correspond
to $\mathcal{N}_a$ and that the edges are ordered so that the edges joining vertices both in $\mathcal{N}_a$
appear last. (There are $\frac{1}{2}|\mathcal{N}_a|(|\mathcal{N}_a|-1)$) such edges, and since $|\mathcal{N}_a| \geq 3$, there are
necessarily at least 3). Quite obviously, for some matrices $A$ and
 $B$ with $n|\mathcal{N}_a|$ columns and with $B$ with at least $\frac{1}{2}|\mathcal{N}_a|(|\mathcal{N}_a|-1)$ rows, there holds
$$
\vspace{-0.2cm}
R=\left[\begin{array}{cc}
R_r & A \\
0 & B
\end{array}\right]
%\vspace{-0.1cm}
$$
Moreover, define $\bar{R}_r$  as a revised reduced rigid matrix,  where the non-zero elements $x_{i1}-x_{l1}$, $x_{i2}-x_{l2}$ in the row of ${R}_r$ which correspond to edges between anchor $l$ and non-anchor $i$  are replaced by $x_{i1}-x_{l1}-\epsilon_{l1}$, $x_{i2}-x_{l2}-\epsilon_{l2}$. 

%  \begin{figure}[ht]
% 	\centering	
% 	\includegraphics[scale=0.15]{demo2.png}\\
%  \vspace{-0.2cm}
% 	\caption{SNL with two non-anchor nodes and four anchor nodes}
% 	\label{fig123}
% 	\vspace{-0.3cm}
% \end{figure}

% For example, for the graph
% of Fig. \ref{fig123}, ${R}_r$ and 
% $\bar{R}_r$ are 
% $$
% R_r=\left[\begin{array}{cccc}
% x_{11}-x_{21} & x_{12}-x_{22} & x_{21}-x_{11} &  x_{22}-x_{21} \\
% x_{11}-x_{31} & x_{12}-x_{32} & 0 & 0  \\
% x_{11}-x_{51} & x_{12}-x_{52}& 0 &0  \\
% 0 & 0 & x_{21}-x_{41} & x_{22}-x_{42} 
% \\
% 0 & 0 & x_{21}-x_{61} & x_{22}-x_{62}
% \\
% \end{array}\right]
% $$
% {\small$$
% \begin{aligned}
%     &\bar{R}_r=\\
%     &\left[\begin{array}{cccc}
% x_{11}\!-\!x_{21} \!\!\!&\!\!\!\!\! x_{12}\!-\!x_{22} \!&\!\!\!\! x_{21}\!-\!x_{11} \!&\!\!\!  x_{22}\!-\!x_{21} \\
% x_{11}\!-\!x_{31}\!-\!\epsilon_{31} \!\!\!&\!\!\! x_{12}\!-\!x_{32}\!-\!\epsilon_{32} \!&\! 0 \!\!\!&\!\!\! \!0  \\
% x_{11}\!-\!x_{51}\!-\!\epsilon_{51} \!\!\!\!&\!\!\! x_{12}\!-\!x_{52}\!-\!\epsilon_{52}\!&\! 0 \!&\!0  \\
% 0 \!&\! 0 \!&\! x_{21}\!-\!x_{41}\!-\!\epsilon_{41} \!\!\!&\!\!\! x_{22}\!-\!x_{42}\!-\!\epsilon_{42} 
% \\
% 0 \!\!\!&\!\!\!\! 0 & x_{21}\!-\!x_{61}\!-\!\epsilon_{61} \!\!\!&\!\!\! x_{22}\!-\!x_{62}\!-\!\epsilon_{62}
% \end{array}\right]
% \end{aligned}
% $$}

% $$
% \bar{R}=\left[\begin{array}{cc}
% R_r & A \\
% 0 & 0
% \end{array}\right]
% $$

Let $\boldsymbol{\rho}$ denote the vector of quantities $\rho_{ij}=\|x_i-x_j\|^2-d_{ij}^{\star}-\mu_{ij}$ and $\rho_{il}=\|x_i-x_l-\epsilon_l\|^2-d_{il}^{\star}-\mu_{il}$ for every 
$(i,j)\in  \mathcal{E}_{ss}$ and
$(i,l)\in  \mathcal{E}_{as}$, which are
with the same ordering as the rows of the revised reduced rigidity matrix $\bar{R}_r$; the entries
depend on the $x_i$, $\mu_{ij}$, $\epsilon_l$ and $\mu_{il}$. Define $\rho_{ij}=0$ ($\rho_{il}=0$) when $(i,j)\notin  \mathcal{E}_{ss}$ 
($(i,l)\notin  \mathcal{E}_{as}$), and let $\Lambda$ denote the
square $N\times N$ matrix with
$$
{\Lambda}\!=\!\left(\begin{array}{ccccc}
-\sum_j \!\rho_{1 j} & \rho_{12} & \rho_{13} & \ldots & \rho_{1 N} \\
% \rho_{12} & -\sum_j \!\rho_{2 j} & \rho_{23} & \ldots & \rho_{2 N} \\
% \rho_{13} & \rho_{23} & -\sum_j\! \rho_{3 j} & \ldots & \rho_{3 K} \\
\vdots & \vdots & & & \vdots \\
\rho_{1 N} & \rho_{2 N} & \rho_{3 N} & \ldots & -\sum_j \!\rho_{j N}
\end{array}\right)
$$

On this basis, we have the following result.
\begin{lemma}\label{l23}
The column vector $\nabla {\Phi}$
whose $(2i-1)_{th}$ and $2i_{th}$ entries are $\frac{\partial {\Phi}}{\partial x_{i1}}$ and $\frac{\partial {\Phi}}{\partial x_{i2}}$ is given by
%\vspace{-0.15cm}
$
\nabla {\Phi} =4\bar{R}_r^T \boldsymbol{\rho}.
$
Further, the Hessian matrix $\nabla^2 \Phi$ is given by 
$$
%\vspace{-0.2cm}
\nabla^2 \Phi=4(2\bar{R}_r^T \bar{R}_r-{\Lambda}\otimes I_{2}).
$$
\end{lemma}

%\noindent\proof
 % \noindent\textbf{Proof sketch}  
  To establish the first claim in Lemma \ref{ee1},   we can apply the implicit function theorem \cite{krantz2002implicit} to the equation $\nabla {\Phi}(\boldsymbol{x}, \boldsymbol{\mu}, \boldsymbol{e}, \boldsymbol{\epsilon})=0$. 
  %The proof of the first claim is given in Appendix \ref{rea}.
  The first claim of the lemma will hold, provided the Jacobian of $\nabla {\Phi}$ is nonsingular at the solution point of $\nabla {\Phi}(\hat{\boldsymbol{x}})=0$ defined by the pair  $(\hat{\boldsymbol{x}}, \boldsymbol{\mu}, \boldsymbol{e}, \boldsymbol{\epsilon})=0$. 

On this basis,
we first show that $\nabla^2 {\Phi}$ is positive definite at the solution point of $\nabla {\Phi}(\boldsymbol{x}^{\star})=0$ defined by the pair    $(\boldsymbol{x}^{\star}, \boldsymbol{\mu}_0, \boldsymbol{e}_0, \boldsymbol{\epsilon}_0)=0$.  
  % $(\boldsymbol{x}^{\star}, \boldsymbol{\mu}_0, \boldsymbol{e}_0, \boldsymbol{\epsilon}_0)$ with $\boldsymbol{\mu}_0=0, \boldsymbol{e}_0=0, \boldsymbol{\epsilon}_0=0$. This proof   can be found in Appendix \ref{rea} and the Jacobian  $\nabla^2 {\Phi}$ is evaluated in Lemma \ref{l23} in Appendix \ref{rea}. 
  %of the positive definiteness of $\nabla^2 {\Phi}$ 
%Hence, $\nabla^2 {\Phi}$ is nonsingular, as required.
%we can prove that  $\nabla^2 {\Phi}$ is positive definite at $(\boldsymbol{x}^{\star}, \boldsymbol{\mu}_0, \boldsymbol{e}_0, \boldsymbol{\epsilon}_0)$.
We can verify that  $\bar{R}_r={R}_r$ has full column rank at $(\boldsymbol{x}^{\star}, \boldsymbol{\mu}_0, \boldsymbol{e}_0, \boldsymbol{\epsilon}_0)$, indicating that $\bar{R}_r^{T}\bar{R}_r$ is positive definite and certainly nonsingular \cite{anderson2010formal}. Also, the matrix $\Lambda$ is zero at this point. Thus,  %$\nabla^2 {\Phi}$ is %positive definite
the statement follows. 

Then, referring to the implicit function theorem \cite{krantz2002implicit},
there exists a neighborhood around $(\boldsymbol{x}^{\star}, \boldsymbol{\mu}_0, \boldsymbol{e}_0, \boldsymbol{\epsilon}_0)$ such that the equation $\nabla {\Phi}({\boldsymbol{x}}, \boldsymbol{\mu}, \boldsymbol{e}, \boldsymbol{\epsilon})=0$ determines ${\boldsymbol{x}}$ as a single-valued function of ${\boldsymbol{x}}=f(\boldsymbol{\mu}, \boldsymbol{e}, \boldsymbol{\epsilon})$ with $\boldsymbol{x}^{\star}=f(\boldsymbol{\mu}_0, \boldsymbol{e}_0, \boldsymbol{\epsilon}_0)$ and $f(\boldsymbol{\mu}, \boldsymbol{e}, \boldsymbol{\epsilon})$ is continuous with continuous partial derivatives with respect to all of its variables. Accordingly, 
$$
\begin{aligned}
\|{\boldsymbol{x}}-\boldsymbol{x}^{\star}\|^2&=(\|f(\boldsymbol{\mu}, \boldsymbol{e}, \boldsymbol{\epsilon})-f(\boldsymbol{\mu}_0, \boldsymbol{e}_0, \boldsymbol{\epsilon}_0)\|)^2\\
&\leq (\|f(\boldsymbol{\mu}, \boldsymbol{e}, \boldsymbol{\epsilon})-f(\boldsymbol{\mu}_0, \boldsymbol{e}, \boldsymbol{\epsilon})\|+\|f(\boldsymbol{\mu}_0, \boldsymbol{e}, \boldsymbol{\epsilon})\\
&\quad-f\!(\boldsymbol{\mu}_0, \!\boldsymbol{e}_0, \!\boldsymbol{\epsilon})\|\!+\!\|f(\boldsymbol{\mu}_0, \boldsymbol{e}_0, \boldsymbol{\epsilon}\!)\!-\!
f\!(\boldsymbol{\mu}_0, \!\boldsymbol{e}_0, \!\boldsymbol{\epsilon}_0)\|\!)^2\\
&\leq a\|\boldsymbol{\mu}\|^2+b\|\boldsymbol{e}\|^2+c\|\boldsymbol{\epsilon}\|^2.
\vspace{-0.25cm}
\end{aligned}
$$
% Also,
% $$
% \begin{aligned}
% {\Phi}(\hat{\boldsymbol{x}})-{\Phi}(\boldsymbol{x}^{\star})
% &\leq l\|\hat{\boldsymbol{x}}-\boldsymbol{x}^{\star}\|^2\\
% &\leq al\|\boldsymbol{\mu}\|^2+bl\|\boldsymbol{e}\|^2+cl\|\boldsymbol{\epsilon}\|^2.
% \end{aligned}
% $$
Thus, for the case with  $\nabla {\Phi}(\hat{\boldsymbol{x}})\!=\!0$ defined by the pair  $(\hat{\boldsymbol{x}}, \boldsymbol{\mu}, \boldsymbol{e}, \boldsymbol{\epsilon})\!=\!0$, if $a\|\boldsymbol{\mu}\|^2+b\|\boldsymbol{e}\|^2+c\|\boldsymbol{\epsilon}\|^2\!< \!\delta$, $\|\hat{\boldsymbol{x}}-\boldsymbol{x}^{\star}\|^2< \delta$. 

  For the second claim,  since the Jacobian  matrix $\nabla^2 {\Phi}$ is positive definite at the point $( \boldsymbol{\mu}_0, \boldsymbol{e}_0, \boldsymbol{\epsilon}_0)$ (shown in Appendix \ref{rea}), the continuity of $\nabla^2 {\Phi}$ in $(\boldsymbol{\mu}, \boldsymbol{e}, \boldsymbol{\epsilon})$ guarantees its positive definiteness for suitably small errors. Thus, the stationary point $\hat{\boldsymbol{x}}$ is a local minimum. Recall the definition of potential game, $\hat{\boldsymbol{x}}$ is also a local NE.
\hfill $\square$

\section{Proof of Theorem \ref{t11}}\label{a4n}

% Let $\mathcal{B}$ denote the ball around $\hat{\boldsymbol{x}}$ defined by $\|\boldsymbol{x}-\boldsymbol{x}^{\star}\|^2< \delta_1$, and let $\mathcal{B}^c$ denote the complementary set $\|\boldsymbol{x}-\boldsymbol{x}^{\star}\|^2\geq \delta_1$. Observe that for all $\boldsymbol{\mu}$, $\boldsymbol{e}$ and $\boldsymbol{\epsilon}$
% with $a\|\boldsymbol{\mu}\|^2+b\|\boldsymbol{e}\|^2+c\|\boldsymbol{\epsilon}\|^2< \delta_1$, Lemma \ref{ee1} guarantees that $\hat{\boldsymbol{x}} \in \mathcal{B}$. Let ${\Phi}_1$
% be defined by
% \vspace{-0.1cm}
% \begin{align}\label{p11}
% {\Phi}_1=&\min _{{x}_i, i \in \mathcal{N}_s, \boldsymbol{x} \in \mathcal{B}^c} 
% %{\Phi}_1(x_1,\dots,x_{{N}})\\&= 
% \sum_{(i, j) \in \mathcal{E}_{ss}} \!\!\!\!{(\left\|{x}_i-{x}_j\right\|^2-d_{i j}^{\star2})^2}\\
% %+  \sum_{l \in \mathcal{N}_a}\!\!\left\|{x}_l-x_l\right\|^2
% &\quad\quad\quad\quad\quad\quad\quad+\!\!\!\!\sum_{(i, l) \in \mathcal{E}_{as}} \!\!\!\!{(\left\|{x}_i-{x}_l^{\star}\right\|^2-d_{i l}^{\star2})^2}. \notag
% %+  \sum_{l \in \mathcal{N}_a}\!\!\left\|{x}_l-x_l\right\|^2
% % &=\sum_{(i, j) \in \mathcal{E}\backslash\mathcal{E}_{aa}} \!\!\!\!{(\left\|{x}_i-{x}_j\right\|^2-(d_{i j}^{\star 2}+\mu_{i j}))^2}+\|\Lambda(\boldsymbol{x}-(\boldsymbol{x}^{\star}+\boldsymbol{\epsilon}))\|^2
% \vspace{-0.35cm}
% \end{align}
Recall the definition of $\Phi_1$ in \eqref{p11}.
The unique localizability property with zero noise, i.e., there exist unique $x_i$ and
$x_j$ such that $\|x_i-x_j\| = d_{ij}^{\star}$, and  $\|x_i-x_l^{\star}\| = d_{ij}^{\star}$, guarantees that $\Phi_1$ is positive.
Also, ${\Phi}_1$ is overbounded by the minimum value of $\Phi$ computed on $\|\boldsymbol{x}-\boldsymbol{x}^{\star}\|^2= \delta_1$.

On the other hand, recall the definition of $\Phi_2$ in \eqref{phi2}.
% Consider also a collection of minimization problems, parameterized by a nonnegative constant $\delta_2$, with variables $\boldsymbol{x}$, $\boldsymbol{\mu}$, $\boldsymbol{e}$ and $\boldsymbol{\epsilon}$:
% \begin{align}\label{phi2}
% \Phi_2=\!\!\!\!\!\!\!\!&\min _{{x}_i, i \in \mathcal{N}_s, \boldsymbol{x} \in \mathcal{B}^c, a\|\boldsymbol{\mu}\|^2+b\|\boldsymbol{e}\|^2+c\|\boldsymbol{\epsilon}\|^2\leq \delta_2} %{\Phi}_2(x_1,\dots,x_{{N}})\\&
% \sum_{(i, j) \in \mathcal{E}_{ss}} \!\!\!\!(\left\|{x}_i-{x}_j\right\|^2\notag\\
% %+  \sum_{l \in \mathcal{N}_a}\!\!\left\|{x}_l-x_l\right\|^2
% &\quad-d_{i j}^{\star2}-\mu_{ij})^2+\!\!\!\!\!\!\!\sum_{(i, l) \in \mathcal{E}_{as}} \!\!\!\!\!\!(\left\|{x}_i\!-\!(x_l^{\star}\!+\!\epsilon_l)\right\|^2\!\!\!-\!(d_{i l}^{\star 2}\!+\!\|\epsilon_l\|^2 \notag\\
% &\quad-\!\!2d_{i l}^{\star}\|\epsilon_l\|\text{cos}(\theta_{il})\!\!+\!\!e_{il}))^2. 
% \vspace{-0.35cm}
% %+\sum_{(i, l) \in \mathcal{E}_{as}} \!\!\!\!{(\left\|{x}_i-{x}_l^{\star}-\epsilon_l\right\|^2-d_{i l}^{\star2}-\mu_{il}-e_{il})^2} \notag
% %+  \sum_{l \in \mathcal{N}_a}\!\!\left\|{x}_l-x_l\right\|^2
% % &=\sum_{(i, j) \in \mathcal{E}\backslash\mathcal{E}_{aa}} \!\!\!\!{(\left\|{x}_i-{x}_j\right\|^2-(d_{i j}^{\star 2}+\mu_{i j}))^2}+\|\Lambda(\boldsymbol{x}-(\boldsymbol{x}^{\star}+\boldsymbol{\epsilon}))\|^2
% \end{align}
There is an infimum as opposed to a minimum used in (\ref{phi2}) because the set over which the extremization is performed is unbounded. Let us now argue that a bounded
set can be used, leading to a minimum result. With $\delta_2$ fixed, we claim that the value of
$\Phi$ in (\ref{potential-fun}) computed on the set $\{(\boldsymbol{x}, \boldsymbol{\mu},\boldsymbol{e},\boldsymbol{\epsilon}): \|\boldsymbol{x}-\boldsymbol{x}^{\star}\|^2 = R \;\text{and} \; a\|\boldsymbol{\mu}\|^2+b\|\boldsymbol{e}\|^2+c\|\boldsymbol{\epsilon}\|^2\leq \delta_2\}$ for large enough
$R$ will go to infinity as $R \rightarrow\infty$.

Thus, the infimum $\Phi_2$ of the index $\Phi$ over the set $\mathcal{B}^c_{\delta_1}\!(\boldsymbol{x}^{\star}) \cap a\|\boldsymbol{\mu}\|^2+b\|\boldsymbol{e}\|^2+c\|\boldsymbol{\epsilon}\|^2\leq \delta_2$ is going
to be attained over the intersection of the set $\mathcal{B}^c_R :=
\mathcal{B}^c_{\delta_1}\!(\boldsymbol{x}^{\star}) \cap \|\boldsymbol{x}-\boldsymbol{x}^{\star}\|^2 \leq R$ for some
suitably large $R$ and $a\|\boldsymbol{\mu}\|^2+b\|\boldsymbol{e}\|^2+c\|\boldsymbol{\epsilon}\|^2\leq \delta_2$, and because this is a bounded and closed set, there
will be at least one point in it achieving the minimum.

Now we choose $\delta_2$. Note that when $\delta_2 =0$, ${\Phi}_2 = {\Phi}_1$.
Because ${\Phi}$ depends continuously on the $\boldsymbol{\mu}$, $\boldsymbol{e}$ and $\boldsymbol{\epsilon}$, ${\Phi}_2$ depends continuously on $\delta_2$ as $\delta_2$ increases from 0. If ${\Phi}_2 \geq \frac{1}{2}{\Phi}_1$ for all $\delta_2 \leq \delta_1$, choose $\delta_2 = \delta_1$. Otherwise,
choose $\delta_2$ so that ${\Phi}_2 = \frac{1}{2}\Phi_1$,  which requires $\delta_2 < \delta_1$.

By the argument above, the global minimum of $\Phi$ over the set $\mathcal{B}^c_{\delta_1}\!(\boldsymbol{x}^{\star}) $ for  fixed $\mu_{ij}$, $e_{il}$, $\epsilon_{l}$
obeying $a\|\boldsymbol{\mu}\|^2+b\|\boldsymbol{e}\|^2+c\|\boldsymbol{\epsilon}\|^2\leq \delta_2$ is at least $\frac{1}{2}{\Phi}_1$. By Lemma \ref{ee1}, there is a single local, and
therefore global, minimum of $\Phi$ in the closure of the set $\mathcal{B}_{\delta_1}\!(\boldsymbol{x}^{\star})$ for any fixed $\mu_{ij}$, $\mu_{il}$, $\epsilon_{l}$ obeying
$a\|\boldsymbol{\mu}\|^2+b\|\boldsymbol{e}\|^2+c\|\boldsymbol{\epsilon}\|^2\leq \delta_1$ and a fortiori $a\|\boldsymbol{\mu}\|^2+b\|\boldsymbol{e}\|^2+c\|\boldsymbol{\epsilon}\|^2\leq \delta_2$. If this minimum is less than $\frac{1}{2}{\Phi}_1$, then it must be the global minimum with no restriction on the set of allowed $\boldsymbol{x}$. Define $\delta$
by 
\vspace{-0.1cm}
$$
\delta=\min\{\delta_2, (\frac{1}{2}{\Phi}_1)\}.
\vspace{-0.15cm}
$$
For any fixed $\boldsymbol{\mu}$, $\boldsymbol{e}$, and $\boldsymbol{\epsilon}$, there holds
\begin{align}\label{pp3}
&\Phi^{\Diamond}=\min _{{x}_i, i \in \mathcal{N}_s} {\Phi}(x_1,\dots,x_{{N}}) \\
&=  \sum_{(i, j) \in \mathcal{E}_{ss}} \!\!\!\!{(\left\|{x}_i-{x}_j\right\|^2-d_{i j}^2)^2}
%+  \sum_{l \in \mathcal{N}_a}\!\!\left\|{x}_l-x_l\right\|^2
+\sum_{(i, l) \in \mathcal{E}_{as}} \!\!\!\!{(\left\|{x}_i-{x}_l\right\|^2-d_{i l}^2)^2} 
%+  \sum_{l \in \mathcal{N}_a}\!\!\left\|{x}_l-x_l\right\|^2
\notag\\
&\leq \sum_{(i, j) \in \mathcal{E}_{ss}} \!\!\!\!{(\left\|{x}_i^{\star}-{x}_j^{\star}\right\|^2-(d_{i j}^{\star 2}+\mu_{i j}))^2}\notag\\
&\quad+\!\!\!\!\!\!\!\sum_{(i, l) \in \mathcal{E}_{as}} \!\!\!\!\!\!(\left\|{x}_i^{\star}\!-\!(x_l^{\star}\!+\!\epsilon_l)\right\|^2\!\!\!-\!(d_{i l}^{\star 2}\!+\!\|\epsilon_l\|^2\!\!-\!\!2d_{i l}^{\star}\|\epsilon_l\|\text{cos}(\theta_{il})\!\!+\!\!e_{il}))^2 \notag \\
&\leq  \! \! \! \! \sum_{(i, j) \in \mathcal{E}_{ss}}\!\!\! \mu_{i j}^2 +\!\!\!\!\! \sum_{(i, j) \in \mathcal{E}_{as}} \!\!\!\!\!(
 \!- \!2({x}_i^{\star}-{x}_j^{\star})^T \epsilon_l \!+  \!2d_{i l}^{\star}\|\epsilon_l\|\text{cos}(\theta_{il})
 \!+ \!e_{i l})^2\notag\\
&\leq  \! \! \!  \sum_{(i, j) \in \mathcal{E}_{ss}} \mu_{i j}^2 + \!\!\!\!\!\sum_{(i, j) \in \mathcal{E}_{as}} (
\|{x}_i^{\star}-{x}_j^{\star}\| \|\epsilon_l\| \!+ \! 2d_{i l}^{\star}\|\epsilon_l\|
 \!+ \!e_{i l})^2\notag\\
&\leq \! \sum_{(i, j) \in \mathcal{E}_{ss}} \! \! \mu_{i j}^2 + \sum_{(i, j) \in \mathcal{E}_{as}} ( 4d_{i l}^{\star}\|\epsilon_l\|
 \!+ \!e_{i l})^2\notag\\
&{\leq  \sum_{(i, j) \in \mathcal{E}_{ss}}\mu_{i j}^2 + \sum_{(i, j) \in \mathcal{E}_{as}}  32d_{i l}^{\star2}\|\epsilon_l\|^2
+2e_{i l}^2}\notag\\
&{\leq  \|\boldsymbol{\mu}\|^2 +  32d_{max}^{\star2}\|\boldsymbol{\epsilon}\|^2
+2\|\boldsymbol{e}\|^2}.\notag
% &=  \sum_{(i, j) \in \mathcal{E}_{ss}} \!\!\!\!{(\left\|{x}_i-{x}_j\right\|^2-d_{i j}^2)^2}
% %+  \sum_{l \in \mathcal{N}_a}\!\!\left\|{x}_l-x_l\right\|^2
% +\sum_{(i, l) \in \mathcal{E}_{as}} \!\!\!\!{(\left\|{x}_i-{x}_l\right\|^2-d_{i l}^2)^2} 
% %+  \sum_{l \in \mathcal{N}_a}\!\!\left\|{x}_l-x_l\right\|^2
% \\
% &\leq \sum_{(i, j) \in \mathcal{E}_{ss}} \!\!\!\!{(\left\|{x}_i^{\star}-{x}_j^{\star}\right\|^2-(d_{i j}^{\star 2}+\mu_{i j}))^2}\\
% &\quad+  \sum_{(i, l) \in \mathcal{E}_{as}} \!\!\!\!{(\left\|{x}_i^{\star}-(x_l^{\star}+\epsilon_l)\right\|^2\!-\!(d_{i l}^{\star 2}+\bar{\mu}_{i l})^2} \\
% &\leq  \sum_{(i, j) \in \mathcal{E}_{ss}} \mu_{i j}^2 + \sum_{(i, j) \in \mathcal{E}_{as}} (
% -2({x}_i^{\star}-{x}_j^{\star})^T \epsilon_l+ \|\epsilon_l\|^2
% -\bar{\mu}_{i l})^2\\
% &\leq  \sum_{(i, j) \in \mathcal{E}_{ss}} \mu_{i j}^2 + \sum_{(i, j) \in \mathcal{E}_{as}} 2(-2({x}_i^{\star}-{x}_j^{\star})^T \epsilon_l+\|\epsilon_l\|^2)^2+2
% \bar{\mu}_{i l}^2\\
% &\leq  \sum_{(i, j) \in \mathcal{E}_{ss}} \mu_{i j}^2 + \sum_{(i, j) \in \mathcal{E}_{as}} 2( 4\|{x}_i^{\star}-{x}_j^{\star})\|^2 \|\epsilon_l\|^2+2\|\epsilon_l\|^4)+2
% \bar{\mu}_{i l}^2\\
% &{\color{red}\leq  \sum_{(i, j) \in \mathcal{E}_{ss}}\mu_{i j}^2 + \sum_{(i, j) \in \mathcal{E}_{as}}  8d_{i l}^{\star2}\|\epsilon_l\|^2+4\|\epsilon_l\|^4
% +2\bar{\mu}_{i l}^2}\\
% &{\color{red}\leq  \|\boldsymbol{\mu}_1\|^2 +  8d_{max}^{\star2}\|\boldsymbol{\epsilon}\|^2+4\|\boldsymbol{\epsilon}\|^2
% +2\|\bar{\boldsymbol{\mu}}_2\|^2}
% %&\leq \delta_1^2+ 32\delta_2^2+2\delta_3^2
% %&< \delta
\end{align}
Take {$a=1$, $b=2$, and $c=32d_{max}^{\star2}$, 
and require $\|\boldsymbol{\mu}\|^2 +  2\|{\boldsymbol{e}}\|^2+32d_{max}^{\star2}\|\boldsymbol{\epsilon}\|^2
 < \delta$.} Then $\Phi^{\Diamond} < (1/2)\Phi_1$. Such a
minimum cannot be achieved in the set $\mathcal{B}^c_{\delta_1}\!(\boldsymbol{x}^{\star})$. However,  it is achieved in the set $\mathcal{B}_{\delta_1}\!(\boldsymbol{x}^{\star})$. Within
this set, there is at most one minimum. Therefore, the minimum in this set is the
global minimum with no restriction on the set of allowed $\boldsymbol{x}$. Thus, the conclusion follows.
\hfill $\square$

\bibliographystyle{IEEEtran}
\bibliography{reference}

\end{document}